\documentclass[12pt,oneside]{amsart}
\usepackage{fullpage}
\usepackage{natbib}
\usepackage{amssymb}
\usepackage{amsmath}
\usepackage{mathtools}
\newtheorem{theorem}{Theorem}[section]
\newtheorem{lemma}[theorem]{Lemma}

\newcommand{\KEYWORDS}[1]{\textit{Keywords:} {#1}\\ }
\newcommand{\HISTORY}[1]{\textit{History:} {#1} }
\newcommand{\ABSTRACT}[1]{\begin{abstract}{#1}\end{abstract}}
\newcommand{\ACKNOWLEDGMENT}[1]{\textbf{Acknowledgement:} {#1} }

\usepackage{graphicx}
\usepackage{subfigure}
\usepackage{color}
\usepackage{enumitem}
\graphicspath{{./figures/}{./}{../}}
\def\betavec{{\boldsymbol\beta}}

\def\B{\mathcal{B}}
\def\X{\mathcal{X}}
\def\R{\mathbb{R}}
\def\betaa{\beta^{Alg}}
\def\K{\mathcal{K}}
\def\U{\mathcal{U}}
\def\P{\mathbb{P}}
\def\E{\mathbb{E}}
\def\H{\mathcal{H}}
\def\RA{\mathcal{R}}
\def\L{\mathcal{L}}
\def\x{\mathbf{x}}
\def\xt{\tilde{\x}}
\def\hist{\{(\x^i,y^i)\}_{i=1}^n}
\def\xUnlab1{{\xt^{1}}}
\def\xUnlabj{{\xt^{j}}}
\def\xUnlabm{{\xt^{m}}}
\def\N{\mathbb{N}}
\def\Nr{\mathcal{N}}
\def\H{\mathcal{H}}
\def\y{\mathbf{y}}
\def\policy{{\boldsymbol{\pi}}}
\def\SigmaBF{\boldsymbol{\Sigma}}
\def\Y{\mathcal{Y}}
\def\u{{\mathbf{u}}}
\def\betaat{\beta^{Alg,\tau}}
\def\ind{\mathbf{1}}
\def\I{\mathcal{I}}
\def\Ialg{I^{\textrm{Alg}}}

\def\F{\mathcal{F}}

\makeatletter
\def\@setcopyright{}
\def\serieslogo@{}
\makeatother

\begin{document}
\author{Theja Tulabandhula and Cynthia Rudin}
\address{MIT}
\email{\{theja,rudin\}@mit.edu}
\title[Robust Optimization using Machine Learning for Uncertainty Sets]{Robust Optimization using Machine Learning for Uncertainty Sets}


\ABSTRACT{
Our goal is to build robust optimization problems for making decisions based on complex data from the past. In robust optimization (RO) generally, the goal is to create a policy for decision-making that is robust to our uncertainty about the future. In particular, we want our policy to best handle the the worst possible situation that could arise, out of an \textit{uncertainty set} of possible situations. Classically, the uncertainty set is simply chosen by the user, or it might be estimated in overly simplistic ways with strong assumptions; whereas in this work, we learn the uncertainty set from data collected in the past. The past data are drawn randomly from an (unknown) possibly complicated high-dimensional distribution. We propose a new uncertainty set design and show how tools from statistical learning theory can be employed to provide probabilistic guarantees on the robustness of the policy.\\
\HISTORY{First version: August 30 2013, current version: \today .}
}


\maketitle


\section{Introduction}
\label{sec:introduction}

In this work, we consider a situation often faced by decision makers: a policy needs to be created for the future that would be a best possible reaction to the worst possible uncertain situation; this is a question of \textit{robust optimization}. In our case, the decision maker does not know what the worst situation might be, and uses complex data to estimate the \textit{uncertainty set}, which is the set of uncertain future situations. 
Here we are interested in answering questions such as: How might we construct a principled uncertainty set from these complex data? 
Can we ensure that with high probability our policy will be robust to whatever the future brings? Can we construct uncertainty sets that are useful for the situation at hand and are not too conservative?

In this paper we address the important setting where detailed data (features) are available to predict each possible future situation. We turn to predictive modeling techniques from machine learning to make predictions, and to define uncertainty sets. 
Models created from finite data are uncertain: given a collection of historical data, there many be many predictive models that appear to be equally good, according to any measure of predictive quality. This was called the \textit{Rashomon effect} by statistician Breiman \citep{breiman-cultures}, and it is this source of uncertainty in learning that we capture while designing uncertainty sets. 

Our concept is possibly best explained through an illustrative example. Consider the minimum variance portfolio allocation problem where our goal is to construct a portfolio of assets. Let us temporarily say that we know exactly what the return for each of the assets in the market will be, and denote $\y \in \Y \subseteq \R^m$ as the vector of these known returns. Let the covariance of the returns be $\SigmaBF$, which is also known in advance. We denote $\policy$ as our choice of portfolio weights. We thus solve the basic decision-making problem:
$$\min_{\policy} \policy^{T}\SigmaBF\policy \;\; \textrm{ s.t. } \policy^{T}\mathbf{1} = 1,\;\; \y^{T}\policy \geq c,$$ 
where  $()^{T}$ is the transpose operator, $c$ is a constant and $\mathbf{1}$ is the vector of all ones. The objective represents the `risk' of the portfolio that we wish to minimize and the two constraints represent that: (a) the sum of portfolio weights should be equal to one, and (b) the return on the portfolio should be lower bounded by an acceptable baseline rate of return denoted by $c$. Now let us consider the more realistic case where the returns $\y$ are not known in advance, and we need to make a decision about portfolio weights $\policy$ under uncertainty (for simplicity of exposition, let us assume that $\SigmaBF$ is known even though in reality we may need to estimate it along with the returns $\y$). If we are able to encode our uncertainty about these forecasted returns using an uncertainty set $\U$, then we can take a robust optimization (RO) approach and solve the following:
\begin{align}
\min_{\policy}  \policy^{T}\SigmaBF\policy \;\;
\textrm{ s.t. }  \policy^{T}\mathbf{1} = 1, \;\; \y^{T}\policy \geq c \; \forall \y \in \U, \nonumber
\end{align}
which gives us a best response to the worst possible outcome $\y$ in uncertainty set $\U$.
The uncertainty set $\U$ can be defined in many ways, and the central goal of this work is how to model $\U$ from complex data from the past. These data take the form of features and labels; for instance in the portfolio allocation problem, the data are $\{(\x^i,\y^i)\}_{i=1}^n$ where an observation $\x^i\in \X \subseteq \R^d$ represents information we could use to predict the returns $\y^i \in \Y$ on past day $i$. These data might include macroeconomic indicators such as interest rates, employment statistics, retail sales and so on, as well as features of the assets themselves. Having complex data like this is very common, but often is not considered carefully within the decision problem. Some of the different ways uncertainty sets can be constructed are:

\noindent
\textbullet \;Using a priori assumptions: We may have \emph{a priori} knowledge about the range of possible future situations. In the portfolio allocation problem, we can assume that we know all possible values of the returns. This knowledge can guide us in constructing the returns uncertainty set $\U$ using interval constraints. That is, $\U := \{\y: \; \forall j \; y_{j}\in [\underline{y}_j,\overline{y}_j]\}$, where we manually select $\underline{y}_j$ and $\overline{y}_j$ for each $j$. Here we ignore the complex past data altogether.

\noindent 
\textbullet \;Using empirical statistics: We could create an uncertainty set using empirical statistics of the data. In the portfolio allocation problem, we might define $\U$ to be the set of all return vectors that are close to return vectors $\y^i$ that have been realized in the past. Or, $\U$ could be the convex hull of past returns vectors. Here we ignore the $\x^i$'s altogether.

\noindent 
\textbullet \;Using linear regression to model complex data: Here, we use the complex past data $\{(\x^i,\y^i)\}_{i=1}^n$, but we make strong (potentially incorrect) assumptions on the probability distribution these data are drawn from. We use these assumptions to define a class of ``good" predictive models $\B$ from $\X\rightarrow \Y$. Then, given a new feature vector $\xt$ (also in $\X$), we use $\B$ to define an ``intermediate'' uncertainty set $\U_{\B}$ of all possible outcomes for each situation $\xt$, and another ``intermediate'' uncertainty set $ \U_{-\B}$ to capture model residuals. Together, these two sets can be used to define $\U$. This is illustrated for the portfolio allocation problem as follows.

We define $\B$ as all linear models $\beta: \X\rightarrow \Y$ that fall in the confidence interval determined using a linear regression fit under the usual normality assumption. We then define $\U_{\B}$ as predicted returns from these ``good" models given a new feature vector $\xt$. Additionally, using past data and normality assumptions, we can define the set of model residuals $ \U_{-\B}$. Finally, $\U_{\B}$ and $\U_{-\B}$ are used to define the set $\U$ in the robust portfolio allocation formulation above.
One should think of $\U_{\B}$ as including all predictions from all models that fit the data reasonably well with respect to the squared loss. And think of $\U_{-\B}$ as the union of prediction intervals around these models. Then, $\U$ is the union of the predictions and prediction intervals from all of the good models. This allows our decision to be robust to future realizations within any prediction interval from any reasonably good model. This approach uses all of the data, but makes strong, possibly untrue assumptions of normality.

\noindent 
\textbullet \; Using machine learning to model complex data, which is the topic of this work: This setting is more general than linear regression and with much weaker assumptions.
Methods that make strong assumptions have limited applicability for modern datasets with thousands of features, and such assumptions may hinder prediction performance. 
In this work, we provide two principled ways to construct set $\U$ using historical data. We will present two methods for each approach. Both of these approaches use tools from statistical learning theory and make minimal assumptions about the data source. In particular:\\

\begin{figure}
 \centering
  \subfigure[Uses optimized set function]{\label{fig:empirical-models-generic-full} \includegraphics[scale=.35]{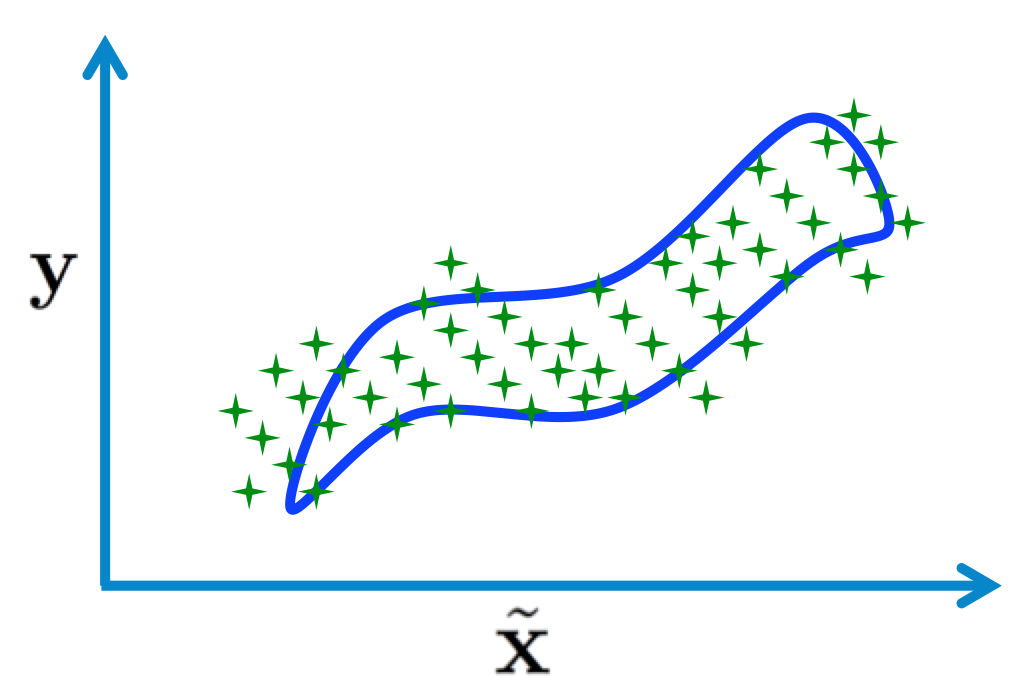}
 }\qquad
 \subfigure[Uses conditional quantile functions]{\label{fig:empirical-models-quantile-full} 
 \includegraphics[scale=.35]{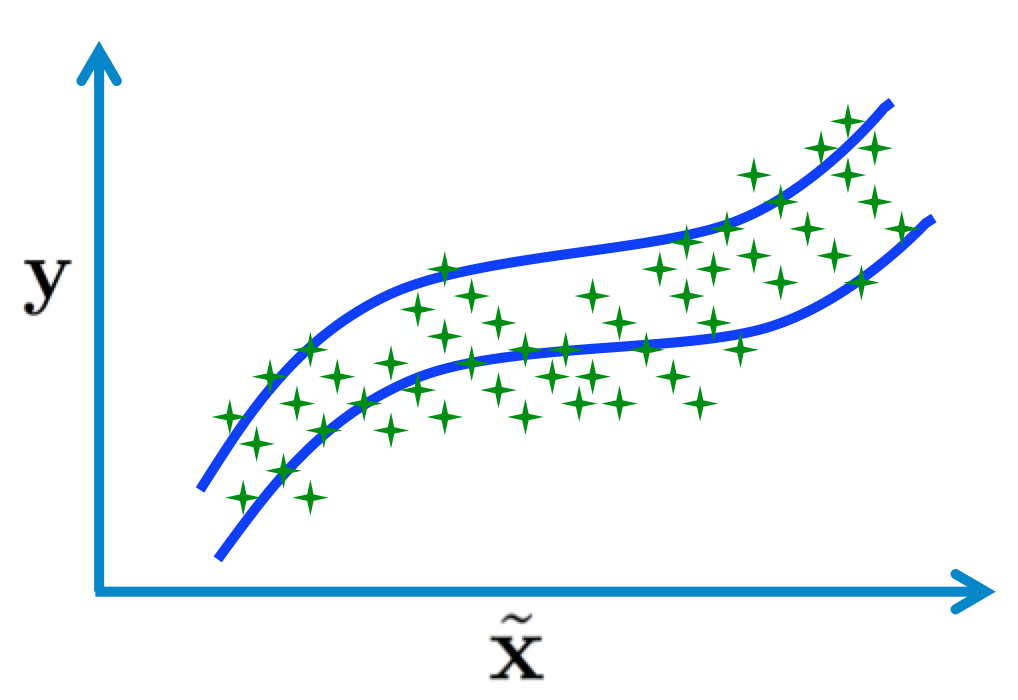}
 }\\
 \subfigure[Uses an intermediate set of ``good'' models]{\label{fig:good-models-mean-full}
 \includegraphics[scale=.35]{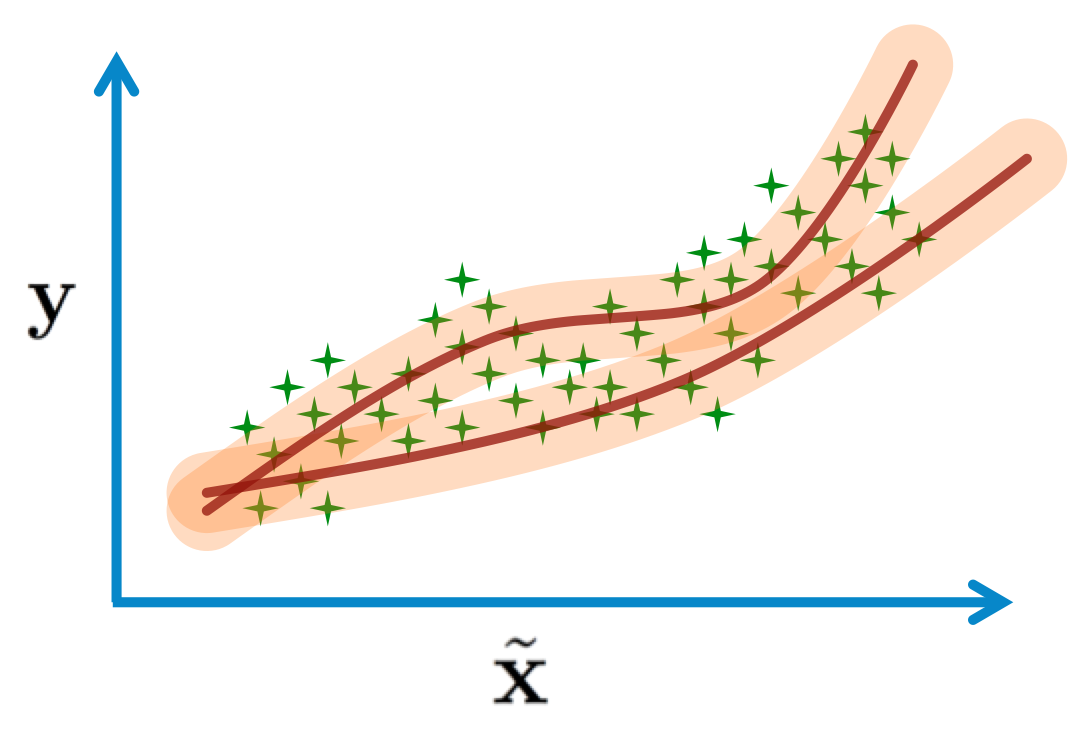}
 }\qquad
 \subfigure[Uses two intermediate sets of ``good'' models]{\label{fig:good-models-quantile-full}
 \includegraphics[scale=.35]{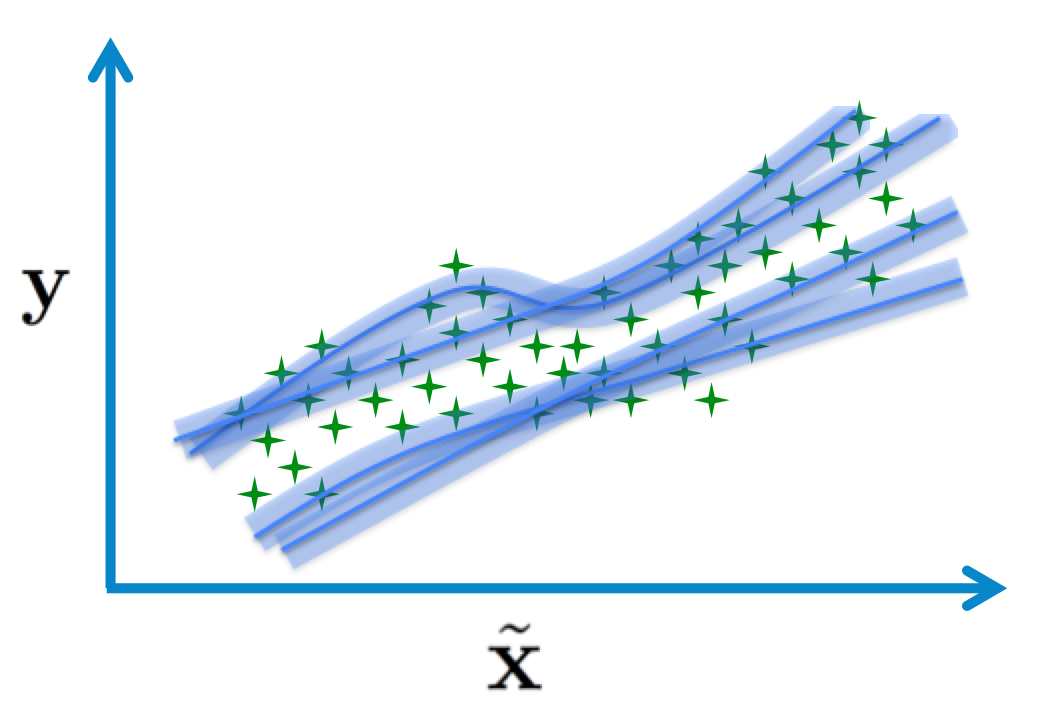}
 }
 \caption{ The empirical data $\{\x^i,\y^i\}_{i=1}^{n}$ is shown along with the boundaries created by the proposed methods in each of the above figures. Evaluation of these boundaries at a given $\xt$ produces an uncertainty set. In (a), a set function is optimized over the sample and its evaluation at every $\xt$ is plotted. In (b), we use optimized conditional quantile models to get the boundaries. In (c), we use an intermediate set of good prediction models and assumptions about model residuals to get the boundaries. In (d), we use two intermediate sets of good conditional quantile models. The lower and upper limits are used to define the boundaries.
 \label{fig:sets-diagram-full}}
 \end{figure}

\noindent\textbf{\textit{(a)}} In the first approach, we optimize prediction models over the data $\{(\x^i,\y^i)\}_{i=1}^n$, and use them to construct uncertainty set $\U$. $\U$ is used within the robust optimization problem to construct $\policy^*$, and Theorem \ref{theorem:empirical-models-generic} provides a guarantee on its robustness; this guarantee is derived using statistical learning theory.
Theorem \ref{theorem:empirical-models-generic} describes the guarantee for a generic class of prediction models and Theorem \ref{theorem:empirical-models-quantile} specializes the guarantee for a specific set of prediction models, namely, the conditional quantile models. Note that in this approach, we do not explicitly construct a set of ``good'' prediction models $\B$ as in the regression approaches discussed in the bullet point above; here $\U$ is defined only from the optimized prediction models and the new feature vector $\xt$. The only assumption made in this approach is that the data are drawn i.i.d from an unknown source distribution. In particular, there is no normality assumption. Let us give examples of how the two methods we propose for this approach would work when $\U$ is constructed from a regression problem (like the portfolio setting discussed earlier):
\begin{itemize}
\item For the first method, for every $\xt$ the uncertainty set $\U$ corresponds to the domain of a indicator function on part of the set $\Y$. It is $1$ on most of the training examples and is $0$ farther away from them. Figure \ref{fig:empirical-models-generic-full} shows an illustration of this.
\item For the second method, we estimate the $95^{\textrm{th}}$ and $5^{\textrm{th}}$ percentiles of $\y$ given $\xt$ and set $\U$ to be all values of $\y \in \Y$ between the two estimates. Figure \ref{fig:empirical-models-quantile-full} illustrates this.
\end{itemize}

\noindent\textbf{\textit{(b)}} In the second approach, we consider the most extreme models within a class of ``good'' models $\B$. The set $\B$ contains all models within a parametric class that have low enough training error. We make only a single assumption: \emph{with high probability, the error due to the `best-in-class' model $\beta^*$ is bounded with a known constant.} Our policies need to be robust to $\beta^*$ that we would choose if we knew the distribution of data. Thus, we make efforts to ensure that the set of good models $\B$ that we will construct contains $\beta^*$. Here, $\B$ and $\U_{\B}$ are chosen in a distribution-independent manner, based on learning theory results. $\U_{-\B}$ is chosen based on our assumption on $\beta^*$. Theorems \ref{theorem:good-models-mean} and \ref{theorem:good-models-quantile} give high probability guarantees on the robust optimal solution obtained using uncertainty set $\U$ constructed in this way. Theorem \ref{theorem:good-models-mean} corresponds to the case where a single prediction model is considered and Theorem \ref{theorem:good-models-quantile} corresponds to the situation where two prediction models (for different quantiles) are considered. These guarantees are qualitatively different from the ones obtained in the first approach. To provide intuition for the two methods proposed for this approach in a regression setting (for instance, as in the portfolio problem):
\begin{itemize}
\item The third method would set $\B$ to be all elements of the hypothesis space (functions on $\X\mapsto \Y$) that have a low least squares loss on the dataset $\{(\x^i,\y^i)\}_{i=1}^n$. These functions estimate the mean of $\y$ given $\xt$. Then we would take an interval above and below each element of $\B$. The union of those intervals would be the uncertainty set $\U$. Figure \ref{fig:good-models-mean-full} illustrates this.
\item The fourth method would set $\B^{0.95}$ to be all models of the $95\textrm{th}$ percentile of $\y$ given $\tilde{x}$ that have low loss. It would set $\B^{0.05}$ to be all models estimating the $5\textrm{th}$ percentile of $\y$ given $\tilde{x}$ that have low loss. We take an interval above and below each estimate provided by $\B^{0.95}$ and $\B^{0.05}$, and take the union of all of these intervals to form $\U$. Note that the fourth method is strictly more conservative than the second method the way we described it. Figure \ref{fig:good-models-quantile-full} illustrates this.
\end{itemize}
 
Being able to define uncertainty sets from predictive models is important: the uncertainty sets can now be specialized to a given new situation $\xt \in \X$, and this is true even if we have never seen $\xt$ before.  For instance, when ordering daily supplies $\y^i$ for an ice cream parlor in Boston, an uncertainty set that depends on the weather might be much smaller than one that does not; planning for too much uncertainty in the weather can be too conservative and very costly: it would not be wise to budget for the largest possible summer sales in the middle of the winter. Though there have been attempts to define uncertainty sets in the linear regression setting \citep{gi03}, ours is the first attempt to tackle the more general setting in a principled way.

Our goals are twofold: (i) We would like to create uncertainty sets for the more general machine learning setting using our proposed approaches (a) and (b) listed above. 
(ii) We would like to compute \textit{sample complexity} values. That is, we want to determine how much data the practitioner needs for a guarantee that their chosen policy will be robust to future realizations. We provide finite sample guarantees on the quality of robustness using learning theory for both proposed approaches. 

Our approaches for constructing uncertainty sets are flexible, intuitive, easy to understand from a practitioner's point of view, and at the same time can bring all the rich theoretical results of learning theory to justify the data-driven methodology. Our uncertainty set designs can handle prediction models for classification, regression, ranking and other supervised learning problems. A main theme of this work is that RO is a new context in which many learning theory results naturally apply and can be directly used.

In Section \ref{sec:formulation}, we formulate our problem and discuss the two approaches (a) and (b) for making decisions under learning uncertainty. In Sections \ref{sec:empirical-functions-general}, \ref{sec:quantiles-without-assumption-b},  \ref{sec:good-models-mean} and \ref{sec:good-models-quantiles}, we use learning theory techniques to justify the proposed uncertainty sets and state our probabilistic guarantees. Section \ref{sec:proofs} provides proofs for these guarantees. Finally, we conclude in Section \ref{sec:conclusion}.

\section{Background Literature}

There are many approaches to decision making under uncertainty when the uncertainty is due to finite data. Robustness is achieved either by taking into the uncertainty in the decision making formulation (as in RO discussed below), or by building robust statistical estimators \citep[see][for applications to portfolio problems]{frost1986empirical,jorion1986bayes}.

In the optimization literature, there has been a continued interest in modeling uncertainty sets for robust optimization (RO) using empirical statistics of data  \citep{delage2010distributionally}, along with (strong) a priori assumptions about the probability distribution generating the parameters of a particular model for the data. \citet{vishal13a} explore a way to specify data-driven uncertainty sets with probabilistic guarantees, where statistical hypothesis testing is used to construct sets. This approach is different from our approach in three important ways: (i) the method is designed for non-complex featureless data, (ii) the goal is totally different: For \cite{vishal13a}, the goal is to minimize the difference between the cost from a policy created using the true distribution and the cost from a policy from the estimated distribution, and (iii) our analysis based on learning theory \citep{vapnik1998statistical} whereas their analysis is based on the theory of hypothesis testing. For us, the objective is to evaluate the feasibility of our policy with respect to a realization of the randomness in the future. The definition of ``robustness" between our work and theirs is thus entirely different.

The closest work to ours is possibly that of \citet{gi03}, who provide a linear-regression-based robust decision making paradigm for portfolio allocation problems, where they assume a multivariate linear regression model for the learning step. A big departure from this approach is that in our work, we are able to design uncertainty sets for a general class of decision making problems while making weak assumptions about the distributional aspects of the historical data. We base our uncertainty set design on regularized empirical risk minimization, which is quite a bit more general than regression. We contrast the sets constructed by \citet{gi03} with our proposed sets in Section \ref{subsec:dependent}.

Our work has the same flavor as \textit{chance constrained programming} \citep{charnes1959chance} and various other \textit{stochastic programming} techniques. 
Both stochastic programming and robust optimization have extensions, for instance, for multi-stage decision making. We focus on single stage optimization. In our previous work \citep{TulabandhulaRu13,TulabandhulaRu14} we considered statistical learning theory bounds also for cases when unlabeled points were available. In that work, we considered prior knowledge about the outcome of an optimization problem that uses the $\tilde{y}^j$s. We showed that this kind of prior knowledge can create better generalization guarantees. Here, instead we study feasibility of the $\tilde{y}^j$s.

\section{Formulation}
\label{sec:formulation}

In Sections \ref{subsec:empirical-function-based} and \ref{subsec:good-models-based} we will describe four ways to construct uncertainty set $\U$ using historical data and solve the corresponding robust optimization problems. The first two methods correspond to approach (a) in the introduction (Section \ref{sec:introduction}), and the last two methods correspond to approach (b). 

Let all the uncertain parameters of the decision problem be denoted by a vector $\u \in \R^m$.
Given a realization of $\u$, let the (basic non-robust) decision making problem be written as:
\begin{align}
\min_\policy \rho(\policy,\u)\;\;\; \textrm{s.t.} \;\; F(\policy,\u) \in \K.\label{eqn:detopt}
\end{align}
Here $\policy \in \Pi \subseteq \R^{d_1}$ is the decision vector and $f:\Pi\times \R^{m} \rightarrow \R$ is the objective function. Function $F:\Pi\times\U \rightarrow \K$ and convex cone $\K \subseteq \R^{d_2}$ describe the constraints of the problem. 

The robust version of the decision problem in Equation (\ref{eqn:detopt}) is thus:
\begin{align}
\min_\policy \max_{\u \in \U} f(\policy,\u)\;\;\; \textrm{s.t.}\;\;  F(\policy,\u) \in \K \textrm{ for all } \u \in \U, \label{eqn:romlopt}
\end{align}
where $\U \subset \R^{m}$ represents the uncertainty set.  In Section \ref{sec:introduction}, the minimum variance portfolio allocation problem is a specific instance of the decision problem in Equation (\ref{eqn:detopt}). The robust portfolio allocation problem is an instantiation of the robust formulation in Equation (\ref{eqn:romlopt}).

To solve Equation (\ref{eqn:romlopt}), we prescribe the following steps:

\noindent\textbf{Step 1}: Construct $\U$ using any of the four methods listed in this section.\\
\noindent\textbf{Step 2}: Obtain a robust solution, using either of the two options below:
\begin{enumerate}[topsep=0pt,leftmargin=*,labelindent=10mm]
\item[Option 1:] 
If $\U$ is a ``nice" set, then there are natural ways  \citep{ben2009robust} to transform it into a relaxed set $\U'$ so that the robust optimization problem can be solved to obtain a robust solution $\pi^*$. For instance, if $U$ can be bounded using a box or an ellipsoid, that box or ellipsoid can be $\U'$. If Equation \ref{eqn:romlopt} is a semi-infinite formulation that can be transformed into a finite formulation, then the finite formulation can be solved.
\item[Option 2:] If $\U$ is not a ``nice'' set, then do the following: sample $L$ elements from $\U$ uniformly. For instance, this can be done using geometric random walks \citep[e.g.,][]{vempala2005geometric}  if $\U$ is convex. Then solve the sampled version of Equation (\ref{eqn:romlopt}) to obtain a robust solution $\policy^*$ \citep[see][]{calafiore2005uncertain} - this method assumes we have an efficient procedure to sample from $\U$.
\end{enumerate}

We focus on \textbf{Step 1}. The goal is to ensure that the true realization of parameter $\u \in \R^m$ belongs to set $\U$ with a high likelihood. Let $\u$ be equal to an $m$-dimensional vector of unknown labels $[\tilde{y}^1 \;\hdots\; \tilde{y}^m]^T$, where each label $\tilde{y}^j \in \Y$ can be predicted
given a corresponding feature vector $\xt^j \in \X$. Thus $m$ labels $\{\tilde{y}^j\}_{j=1}^{m}$, which can be forecasted from $\{\xt^j\}_{j=1}^{m}$, feed into the decision problem of Equation (\ref{eqn:romlopt}).

In both approaches we propose, we will define $\U$ to be a product of $m$ sets, each one constructed such that it contains the corresponding unknown true realization $\tilde{y}^j$ with high probability. Set $\U$ will be a function of training data sample $S = \{\x^i,y^i\}_{i=1}^{n}$ and the current feature vectors $\{\xt^j\}_{j=1}^{m}$.

\subsection{Direct use of empirically optimal prediction models}
\label{subsec:empirical-function-based}
In this approach, we use empirically optimal prediction models directly. We start by discussing a very general form of prediction model, then discuss quantile regression.

\begin{figure}
 \centering
 \subfigure[Using an optimized set function]{\label{fig:empirical-models-generic-slice} 
 \includegraphics[scale=.25]{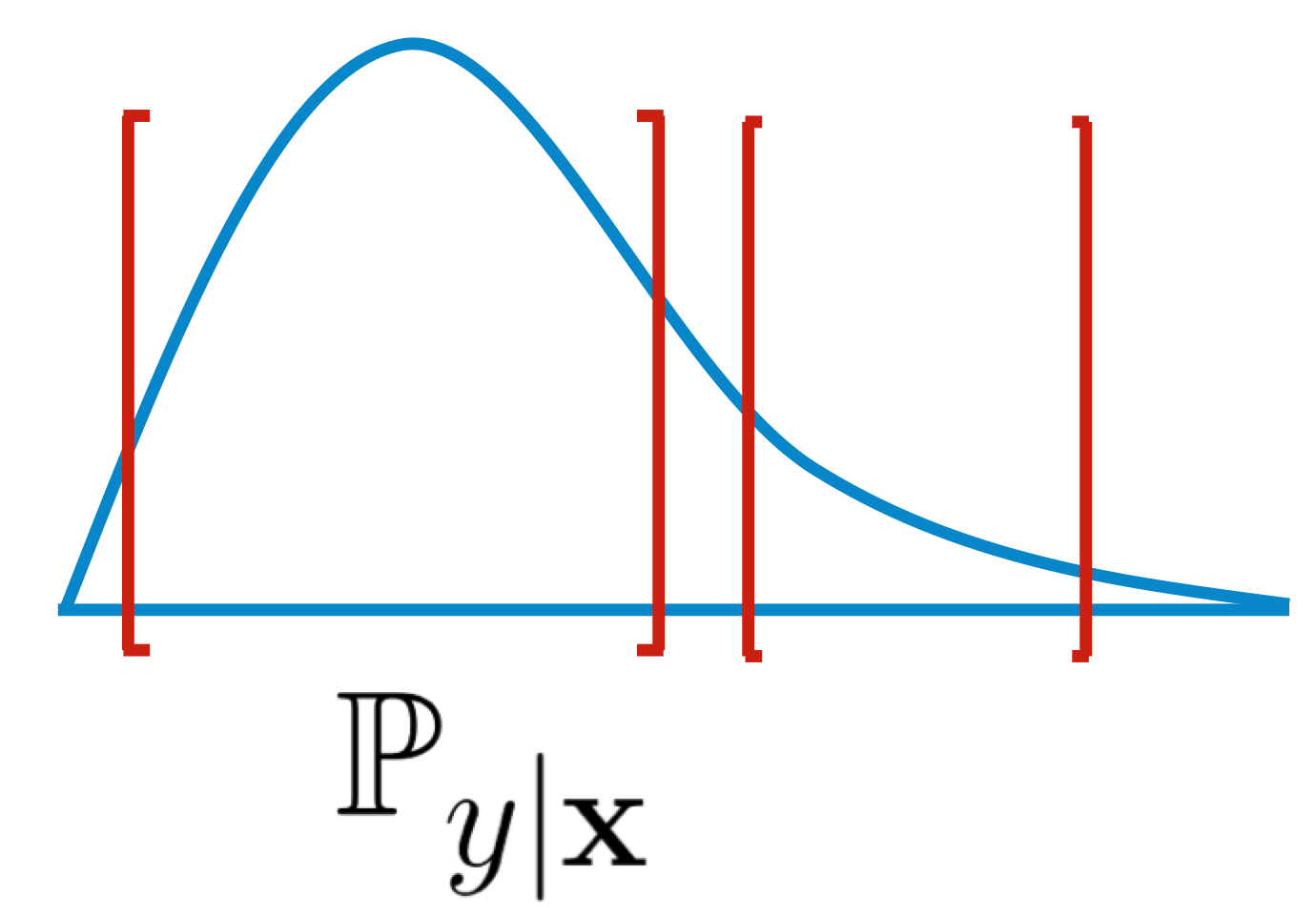}
 }
 \subfigure[Using optimized conditional quantile functions]{\label{fig:empirical-models-quantile-slice} 
 \includegraphics[scale=.25]{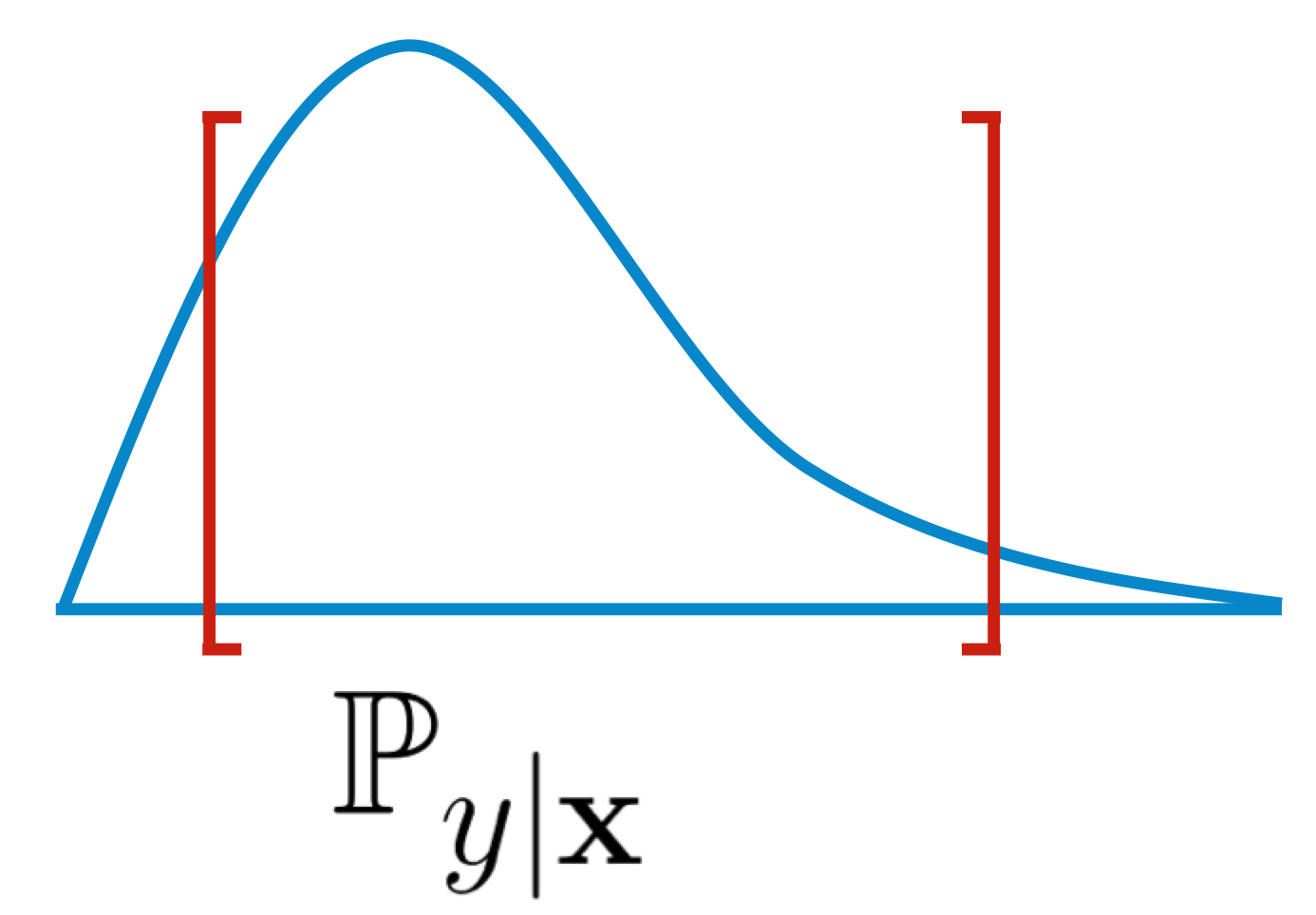}
 }\\
 \subfigure[Using a single intermediate set of ``good'' models]{\label{fig:good-models-mean-slice}
 \includegraphics[scale=.25]{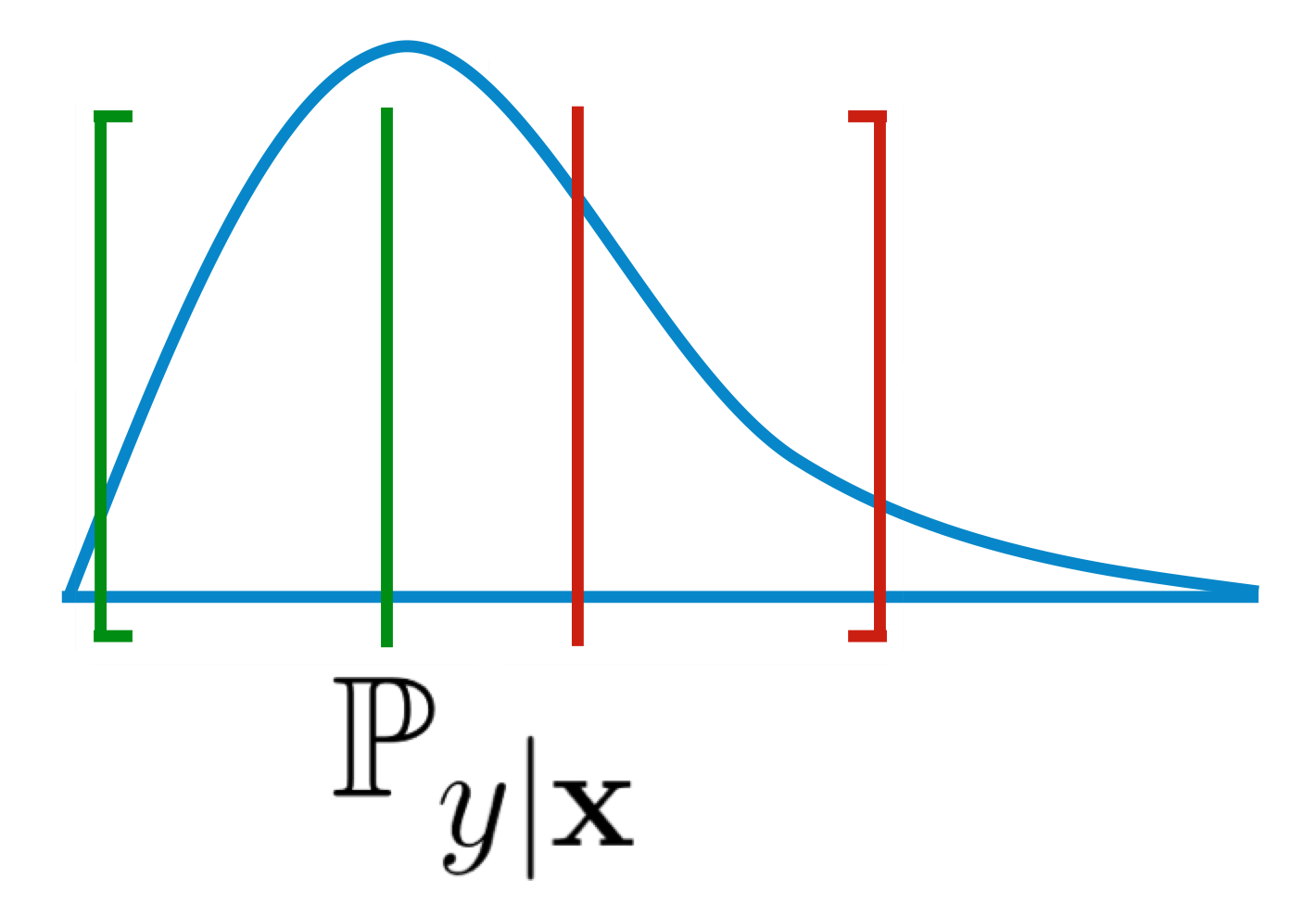}
 }\qquad
 \subfigure[Using two intermediate sets of ``good'' models]{\label{fig:good-models-quantile-slice}
 \includegraphics[scale=.25]{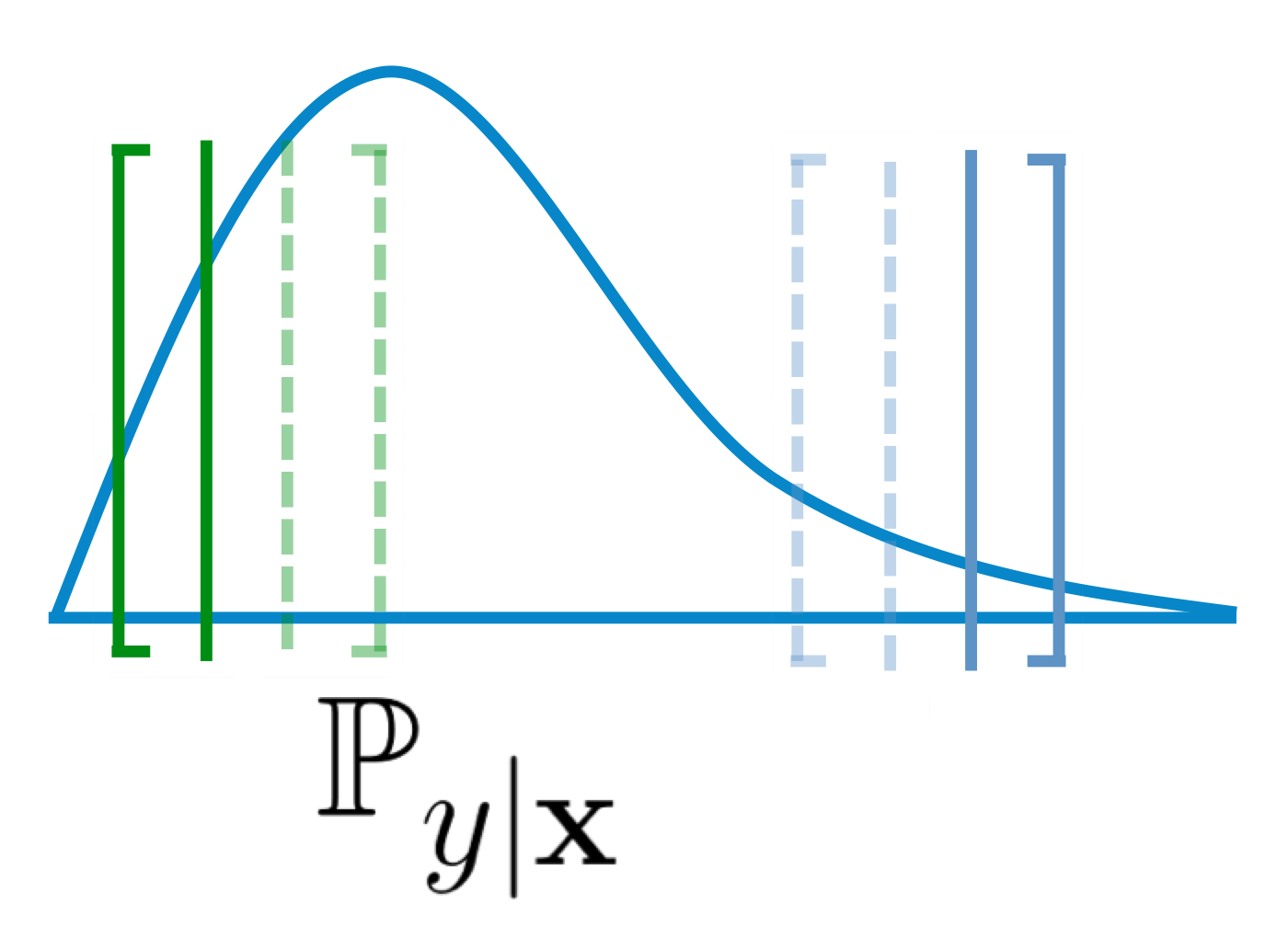}
 }
 \caption{The conditional distribution of $y$ given $\x$ is shown along with the proposed uncertainty sets in each of the above figures. In (a), we use an optimized set function to directly define the subset of $\R$ that contains $y$ with high probability. In (b), we use optimized conditional quantile models (the ones achieving the lowest training error) to directly define the set which contains the random variable $y$ with high probability. In (c), we use an intermediate set of good prediction models to create $\U_{\B}$ and then enlarge the interval using set $\U_{-\B}$. In (d), we use two intermediate sets of good conditional quantile models and enlarge the corresponding intervals. The lower and upper limits of the two sets are then used to define $\U$.
 \label{fig:sets-diagram-slice}}
 \end{figure}

\noindent\textbf{General prediction models}:

Let $\x \in \X \subset \R^d$ represent a feature vector and $y \in \Y \subseteq \R$ represent a label. Consider a class of set functions $I \in \I$, where $I:\X\rightarrow \mathfrak M_\R$, where $\mathfrak M_\R$ is the set of all measurable sets of $\R$.  Let us say that we have a procedure that picks a function $\Ialg$ so that most of the labels of the training examples obey $y^i \in \Ialg(\x^i),\; i=1,...,n$. As long as $\Ialg$ belongs to a set of ``simple'' functions, we have a guarantee on how well $\Ialg$ will generalize to new observations. 
Specifically, consider the following empirical risk minimization procedure:
\begin{align}
\min_{I \in \I}\frac{1}{n}\sum_{i=1}^{n}\ind[y^{i} \notin I(\x^i)], 
\label{eqn:set-erm}
\end{align}
where $\ind[\cdot]$ is the indicator function. Let an optimal solution to the above problem be $\Ialg$. Then, define the uncertainty set $\U$ as: 
\begin{align}
\U  = \Pi_{j=1}^{m}\Ialg(\xt^{j}),
\label{eqn:empirical-models-general-set}
\end{align}
where $\U$ is a product of $m$ measurable sets. Figure \ref{fig:empirical-models-generic-slice} illustrates this construction in one dimension. Given this construction, \textbf{Step 1} of the workflow we described can be summarized as:
\begin{enumerate}
\item[(a)] Solve Equation (\ref{eqn:set-erm}) to obtain a set function $\Ialg$ that depends on sample $S$. 
\item[(b)] Define $\U$ according to Equation (\ref{eqn:empirical-models-general-set}) using new observations $\{\xt^j\}_{j=1}^{m}$.
\end{enumerate}

The above setting is quite general. In particular, since the range of function $\Ialg$ is $\mathfrak M_\R$, we can capture sets that are arbitrarily more complicated than simple intervals. For instance, if $\P_{y^j|\xt^j}$ is bimodal, then for certain values of $\xt^j$, $\Ialg(\xt^j)$ can be the union of two disjoint intervals.

We remark that one can also approximate the source distribution $\P_{\x,y}$ using an empirical distribution $\hat{\P}_{\x,y}$ (there are many parametric and non-parametric ways to do this) and then construct set $\U$ using marginal distributions $\{\hat{P}_{y^j|\xt^j}\}_{j=1}^{m}$. This would be slightly different than the approach described above in that it would require density estimation, which may itself be a hard problem. In the above method and the other methods below that we propose, we focus on estimating  functionals of the conditional distributions $\{\hat{P}_{y^j|\xt^j}\}$ directly.

\noindent\textbf{Conditional quantile models}:

In this method, we specialize the generic function class $\I$ to the class of set functions defined using conditional quantile models. We will estimate an upper quantile of $\tilde{y}$ for each $\xt$, and a lower quantile of $\tilde{y}$ for each $\xt$. The uncertainty set will be the interval between the two quantile estimates. This method is applicable when our prediction task is a regression problem. 

When $y \sim \P_{y}$, the $\tau^{th}$ quantile of $y$, denoted by $\mu^\tau$, is defined as
$
\mu^{\tau} := \inf\{\mu :\; \P_{y}(y \leq \mu) = \tau\}.
$
Here $\tau$ can vary between $0$ and $1$. In the special case when $\tau$ is set to $0.5$, this defines the median. Similarly, when $(\x,y) \sim \P_{\x,y}$, the conditional quantile $\mu^{\tau}$ can be defined as a function from $\X$ to $\Y$,
$
\mu^{\tau}(\x) := \inf\{\mu : \P_{y|\x}(y \leq \mu) = \tau\}.
$

In our setting, $\tilde{y}^j$ conditioned on $\xUnlabj$  is distributed according to $\P_{\tilde{y}^j|\xUnlabj}$. Thus, given a value of $\tau \in [0,1]$, $\P_{\tilde{y}^j|\x = \xUnlabj}(\tilde{y}^j \leq \mu^{\tau}(\xUnlabj)) = \tau$ where $\mu^{\tau}(\x)$ is the conditional quantile defined earlier. Our method picks two values of $\tau$, $\delta_p \leq \delta_q$ such that:
\begin{align*}
\P_{\tilde{y}^j|\xUnlabj}(\tilde{y}^j \leq \mu^{\delta_p}(\xUnlabj)) = \delta_p, \;\textrm{ and }\;
\P_{\tilde{y}^j|\xUnlabj}(\tilde{y}^j \leq \mu^{\delta_q}(\xUnlabj)) = \delta_q.
\end{align*}
For example, a typical value for the pair $(\delta_p,\delta_q)$ can be $(0.05,0.95)$ which makes  $\mu^{\delta_p}(\xUnlabj)$ correspond to the $5$\% conditional quantile and $\mu^{\delta_q}(\xUnlabj)$ correspond to the $95$\% conditional quantile. Given these two conditional quantiles, we have:
\begin{align*}
\P_{\tilde{y}^j |\xUnlabj}(\mu^{\delta_p}(\xUnlabj) < \tilde{y}^j  \leq \mu^{\delta_q}(\xUnlabj)) = \delta_q - \delta_p.
\end{align*}
Thus, the unknown future realization of $\tilde{y}^j$ belongs to the interval  $[\mu^{\delta_p}(\xUnlabj),\mu^{\delta_q}(\xUnlabj)]$ with high probability if $\delta_p$ and $\delta_q$ are chosen appropriately. If we knew the true conditional quantiles (which we do not), we could define the uncertainty set $\U$ as
$\U = \Pi_{j=1}^{m}[\mu^{\delta_p}(\xUnlabj),\mu^{\delta_q}(\xUnlabj)].$
We will circumvent this issue by using sample $S = \hist$ and quantile regression to obtain empirical quantile functions. 

Quantile regression can be seen as an empirical risk minimization algorithm where the loss function is defined appropriately to obtain a conditional quantile function. That is, we aim to obtain an estimator function $\beta(\x)$  of the true conditional quantile function $\mu^{\tau}(\x)$ given a predefined quantile parameter $\tau$.
In particular, the \textit{pinball} loss (or newsvendor loss) function defined below is used.
\begin{align*}
l^{\tau}(\beta(\x),y) = \begin{cases} \tau\cdot(y-\beta(\x)) & \textrm{if } y - \beta(x) \geq 0,\\
(\tau - 1)\cdot(y-\beta(\x)) & \textrm{otherwise.}
\end{cases}
\end{align*}
Let $l_{\P}^{\tau}(\beta) = \E_{\x,y}[l^{\tau}(\beta(\x),y)]$. It can be shown \citep{koenker2005quantile,takeuchi2006nonparametric} under some regularity conditions that the true conditional quantile function $\mu^{\tau}(\x)$ is the minimizer of $l_{\P}^{\tau}(\beta)$ when minimized over all measurable functions. There are several works that consider linear and nonparametric quantile estimates using this loss function \citep{takeuchi2006nonparametric,gahyi14}. In our setting, we will let $\B_0$ be our hypothesis class that we want to pick conditional quantile functions from. 

Let the empirical risk minimization procedure using the pinball loss output a conditional quantile model $\betaat$ when given the historical sample $S = \hist$ of size $n$ and a parameter $\tau$. That is, let $l_{S}^{\tau}(\beta) = \frac{1}{n}\sum_{i=1}^n l^{\tau}(\beta(\x^i),y^i)$ and  $\betaat \in \arg\min_{\beta \in \B_0}l_{S}^{\tau}(\beta)$. The following definition of $\U$ uses two empirical conditional quantile functions with $\tau = \delta_p$ and $\tau = \delta_q$ respectively:
\begin{align}
\U = 
\Pi_{j=1}^{m}\left[\min\left(\beta^{\textrm{Alg},\delta_p}(\xUnlabj),\beta^{\textrm{Alg},\delta_q}(\xUnlabj)\right),\max\left(\beta^{\textrm{Alg},\delta_p}(\xUnlabj),\beta^{\textrm{Alg},\delta_q}(\xUnlabj)\right)\right].
 \label{eqn:empirical-models-quantile-set}
\end{align}
Here $\U$ is again a product of $m$ intervals, each one constructed so that it contains the unknown $\tilde{y}^j$ with high probability (which we prove later). Figure \ref{fig:empirical-models-quantile-slice} illustrates this construction in one dimension. Thus, for \textbf{Step 1}, we do the following:
\begin{enumerate}
\item Compute $\beta^{\textrm{Alg},\delta_p}$ and $\beta^{\textrm{Alg},\delta_q}$ using quantile regression.
\item Set $\U$ according to Equation (\ref{eqn:empirical-models-quantile-set}).
\end{enumerate}

\subsection{Uncertainty set using an intermediate set of ``good'' prediction models}
\label{subsec:good-models-based}

In this approach, we use optimized prediction models to define an intermediate set of ``good'' prediction models, which is then used to define $\U$. This approach aims to capture uncertainty in the modeling procedure explicitly: rather than using one predictive model, we use predictions from all models that we consider to be ``good'' with respect to our training data.

\noindent\textbf{Using a single set of ``good'' prediction models}:

Let $\beta: \X \mapsto \Y$ be a prediction model in the hypothesis class $\B_0$. For instance, $\B_0$ can be the set of linear predictors $\B_0=\{x \mapsto \beta^Tx: \|\beta\| \leq B_b\}$. Let $l(\beta(\x),y)$ denote the loss function. For example, $(\beta(\x)-y)^2$ is the least squares loss and $[1-\beta(\x)y]_+$ is the hinge loss used in Support Vector Machines. For any given model, let $l_{\P}(\beta) = \E_{\x,y}[l(\beta(\x),y)]$ where the expectation is with respect to the unknown distribution $\P_{\x,y}$. Let $\beta^* \in \arg\min_{\beta \in \B_0} l_{\P}(\beta)$ be defined as the `best-in-class' model with respect to our class $\B_0$. Note that we cannot calculate $\beta^*$ as we do not have the distribution.

Our set construction method takes into account two things: (i) how the solution $\betaa$ of empirical risk minimization compares with $\beta^*$ (coming from statistical learning theory), and (ii) how much of the mass of $\P_{\x,y}$ concentrates around $\beta^*(\x)$ (coming from \textbf{Assumption A} described below). 

It is always true that there exists a set $E$ and a scalar $\delta_e \geq 0$ such that:
\begin{align}
\P_{\x,y}\left(\x,y: |y-\beta^*(\x)| \in E \right) \geq 1-\delta_e,
\label{eqn:assumptionA}
\end{align}
where $E \subseteq \Y$. This is trivially satisfied if $E = \Y$. In this case, $\delta_e$ can be set to $0$. Ideally, we know of a pair $(E,\delta_e)$ where $\delta_e$ is still small and where $E$ is not too large; if $E$ were very large, the uncertainty set would be too conservative.  We formalize the assumption that we will use to define $\U$ as follows:

\noindent\textbf{Assumption A:} We know a pair $(E,\delta_e)$ such that Equation (\ref{eqn:assumptionA}) holds.

We can intuitively think of decomposing $\u$ in Equation (\ref{eqn:romlopt}) to capture model uncertainty and residual uncertainty as follows. Let $\u_{\beta}$ be the part of $\u$ that is derived from a statistical model $\beta$. Thus, given $\{\xUnlabj\}_{j=1}^{m}$, $\u_{\beta} := [\beta(\xUnlab1)\;\cdots\beta(\xUnlabm)]^T$. Let the remaining part of $\u$, denoted by $\u_{-\beta}$, be equal to a vector of corresponding model residuals. Thus, $\u = \u_{\beta} + \u_{-\beta}$.

Let $\B$ represent a set of ``good'' prediction models. Let $\U$ be equal to $\U_{\B} + \U_{-\B}$ such that $\u_{\beta} \in \U_{\B}$ and $\u_{-\beta} \in \U_{-\B}$. Here, $\U_{\B}$ corresponds to $\B$ in the following way: $\U_{\B} := \{\u_{\beta}: \beta \in \B\}$. On the other hand, $\U_{-\B}$ corresponds to a set that captures the support of most model residuals. Formally,
\begin{align}
\U = \Pi_{j=1}^{m}\left[\inf\{\beta(\xt^j): \beta \in \B\} -E, \sup\{\beta(\xt^j): \beta \in \B\} + E\right]. \label{eqn:good-models-mean-set}
\end{align}

An illustration in one dimension, when the set of ``good'' models has two members, is shown in Figure \ref{fig:good-models-mean-slice}.
If we know the `best-in-class' model $\beta^*$, then $\U_{\B}$ can be a singleton set just containing $\beta^*$. Since we do not know $\beta^*$, we adapt \textbf{Step 1} of the general recipe to construct $\U$ using $\U_{\B}$ and $\U_{-\B}$ as follows:
\begin{enumerate}
\item[(a)] Define $\B$ using $S = \hist$.  Our sets will be of the form (discussed further in Section \ref{sec:good-models-mean}) 
$\B = \{\beta: g(\beta) \leq g(\betaa) + c\},$
where $g$ is some function, $\betaa$ is a specific model and $c$ is a parameter. These quantities will depend on the learning algorithm and $\hist$.
\item[(b)] Define $\U_\B$ and $\U_{-\B}$: Recall that $\U_{\B} := \{\u_{\beta}: \beta \in \B\}$ where $\u_{\beta} = [\beta(\xUnlab1)\;\cdots\beta(\xUnlabm)]^T$. $\U_{-\B}$ is defined using assumption \textbf{Assumption A} such that it captures the support of the model error residuals (more details are in Section \ref{sec:good-models-mean}). $\U$ is then $\U_{\B} + \U_{-\B}$.
\end{enumerate}

The quality of the robust solution of Equation (\ref{eqn:romlopt}) depends on the set $E$. For a less conservative solution, we want set $E$ to be as small as possible. The probabilistic guarantee on the robust solution that we derive in Section \ref{sec:good-models-mean} depends on $\delta_e$. For a better guarantee, we need $\delta_e$ to be as close as possible to $0$. If our model class $\B_0$ is very complex and able to closely capture most $y$ values, this could reduce the size of set $E$. 

Note that if $\B$ does not contain good models, $\U_{-\B}$ will necessarily be large, our bound on robustness will be loose, and the robust solution thus obtained will be too conservative.

\noindent\textbf{Using two sets of ``good'' prediction models}:

When our prediction problem is a regression task,  we can make a different (and often weaker) assumption than \textbf{Assumption A} using quantile regression. We will construct uncertainty set $\U$ in a different way. 

Recall the definition for the conditional quantile function $\mu^{\tau}(\x)$ and the empirical procedure to estimate it, outlined in Section \ref{subsec:empirical-function-based}. Let $\beta^{\tau,*} \in \arg\min_{\beta \in \B_0} l_{\P}^{\tau}(\beta)$ be the `best-in-class' conditional quantile function for any given $\tau$.  
It is always true that there exists a set $E^{\tau} \subseteq \Y$ and a scalar $\delta_e^{\tau} \geq 0$ such that:
\begin{align}
\P_{\x}(\x: |\mu^{\tau}(\x) - \beta^{\tau,*}(\x)| \in E^{\tau}) \geq 1 - \delta_e^{\tau}.
\label{eqn:assumptionB}
\end{align}
The way we will construct $\U$ below will be such that the quality of the robust solution $\policy^*$ of Equation (\ref{eqn:romlopt}) depends on the set $E^{\tau}$. For a less conservative solution, we want $E^{\tau}$ to be as small as possible. The probabilistic guarantee that we derive in Section \ref{sec:good-models-quantiles} on $\policy^*$ will depend on $\delta_e^{\tau}$. For a better guarantee, $\delta_e^{\tau}$ needs to be as close as possible to $0$. If $\B_0$ is sufficiently rich, $E^{\tau}$ can be small or even empty (which is the case when $\mu^{\tau} \in \B_0$). Thus, similar to \textbf{Assumption A} in Section \ref{subsec:good-models-based}, we make the following assumption:\\
\noindent\textbf{Assumption B}: Given a value of $\tau$, we know a pair $(E^{\tau},\delta_e^{\tau})$ such that Equation (\ref{eqn:assumptionB}) holds.

Let $\B^{\delta_p}$ be the set of ``good'' conditional quantile functions when $\tau = \delta_p$ and let $\B^{\delta_q}$ be the set of ``good'' conditional quantile functions when $\tau = \delta_q$. By ``good'' we mean that all these quantile functions have their quantile estimation performance close to the best we can obtain from $\B_0$ using quantile regression. A precise definition for $\B^{\tau}$ will be given in Equation (\ref{eqn:quantileset}) below. We can then construct $\U$ as:
\begin{align}
\U = \Pi_{j=1}^{m}\Big[&\inf\{\beta(\xUnlabj): \beta \in \B^{\delta_p}\cup\B^{\delta_q}\}-\sup E^{\delta_p}\cup E^{\delta_q}, \nonumber \\
 &\sup\{\beta(\xUnlabj): \beta \in \B^{\delta_p}\cup\B^{\delta_q}\}+\sup E^{\delta_p}\cup E^{\delta_q}\Big].
 \label{eqn:good-models-quantile-set}
\end{align}

The definition of the $j^{\textrm{th}}$ interval involves two sets. The first set, $\{\beta(\xUnlabj): \beta \in \B^{\delta_p}\cup\B^{\delta_q}\}$ contains all the predictions by models in both $\B^{\delta_q}$ and $\B^{\delta_q}$ on the feature vector $\xUnlabj$. The second set, $E^{\delta_p}\cup E^{\delta_q}$, contains all deviations between the true conditional quantiles and the `best-in-class' conditional quantiles at both values of $\tau$. Thus, the smallest value of the predicted $\delta_p$ conditional quantiles and $\delta_q$ conditional quantiles in the first set, in conjunction with the largest deviation captured by the second set, is used to define the lower limit of the interval. The upper limit of the interval is defined in a similar way by taking the largest predicted quantile from the first set and adding the largest deviation captured by the second set. An illustration in one dimension, when each of the two sets of ``good'' models has two members, is shown in Figure \ref{fig:good-models-quantile-slice}.

Given $\U$, we can solve Equation (\ref{eqn:romlopt}) for $\policy^*$ using \textbf{Step 2}. The following is a summary of the way to construct $\U$ in \textbf{Step 1}:
\begin{enumerate}
\item[(a)] Define $\B^{\delta_p}$ and $\B^{\delta_q}$ using $\hist$. We will propose procedures for designing $\B^{\delta_p}$ and $\B^{\delta_q}$ using learning theory results in Section \ref{sec:good-models-quantiles}. Our sets will be of the form
$\B^{\tau} = \{\beta: g(\beta) \leq g(\betaat) + c\}$
for $\tau = \delta_p,\delta_q$, where $g$ is some function, $\betaat$ is a specific conditional quantile model depending on $\tau$ and $c$ is a parameter. These quantities will depend on the pinball loss function and $\hist$.
\item[(b)] Define $\U$:  Using the above sets and the property of quantile error residuals as in Equation (\ref{eqn:assumptionB}), and by \textbf{Assumption B}, we can construct $\U$ as shown in Equation (\ref{eqn:good-models-quantile-set}).
\end{enumerate}

In the next few sections, we provide probabilistic guarantees on the feasibility of the robust optimal solutions obtained by using uncertainty sets from each of the four methods we proposed.

\section{Robustness guarantee using general prediction functions}
\label{sec:empirical-functions-general}

Consider the setting described in Section \ref{sec:formulation}, where we have a class of general set functions $\I$. Let $S := \hist$ be the training data which are independent and identically distributed. Let algorithm $A$ represent a generic learning procedure. That is, it takes $S$ as an input and outputs $\Ialg$. Since $\Ialg$ is a function of sample $S$, we will show that the unknown $\tilde{y}^j$ belong to the interval $\Ialg(\xt^{j})$ with high probability over $S$ as long as the set of functions $\I$ from which $\Ialg$ is picked is ``simple''. Note that we do not assume anything about the source distribution.

In order to state our result, we will define the following quantity known as the empirical Rademacher average \citep{bartlett2002localized}. For a set $\F$ of functions, the \textit{empirical Rademacher average} is defined with respect to a given random sample $S' = \{z^i\}_{i=1}^{n}$ as
$
\RA_{S'}(\F)  = \E_{\sigma^1,...,\sigma^n}\left[\frac{1}{n}\sup_{f \in \F}\left|\sum_{i=1}^{n} \sigma^i f(z^i)\right|\right]
$
where for each $i=1,..,n$, $\sigma^i = \pm 1$ with equal probability. The \textit{Rademacher average} is defined to be the expectation of the empirical Rademacher average over the random sample $S$: $\RA(\H) = \E_{z^1,...,z^n}\left[ \RA_{S}(\H)\right]$. The interpretation of the Rademacher average is that it measures the ability of function class $\F$ to fit noise, coming from the random $\sigma_i's$. If the function class can fit noise well, it is a highly complex class. The Rademacher average is one of many ways to measure the richness of a function class, including covering numbers, fat-shattering dimensions \citep{bartlett1996fat} and the Vapnik-Chervonenkis dimension \citep{vapnik1998statistical}. 

\begin{theorem} If $\U$ is defined as in Equation (\ref{eqn:empirical-models-general-set}), then with probability at least $1-\delta$ over training sample S, we have robustness guarantee
\begin{align*}
\P_{\{\xUnlabj,\tilde{y}^j\}_{j=1}^{m}}\Big(F(\pi^*,[\tilde{y}^1...\tilde{y}^m]^T) \in \K\Big) \geq 
  \left(\left[1 - \frac{1}{n}\sum_{i=1}^{n}\ind[y^i \notin \Ialg(\x^i)] - 2\RA(l\circ \mathcal{I}) - \sqrt{\frac{\log\frac{1}{\delta}}{2n}}\right]_{+}\right)^m,
\end{align*}
where $\epsilon > 0$ is a pre-determined constant, and $\Big[\cdot \Big]_{+}$ is shorthand for $\max(0,\cdot)$.
\label{theorem:empirical-models-generic}
\end{theorem}

The result is a lower bound on the probability of infeasibility. This bound depends on the performance of the data dependent set function $\Ialg$. If $\Ialg$ is such that its performance, measured in terms of $\frac{1}{n}\sum_{i=1}^{n}\ind[y^i \notin \Ialg(\x^i)]$ is good (i.e., lower in value), then the right hand side of the inequality increases, resulting in a higher chance of feasibility. This probability of feasibility also depends on the number of estimates $m$ that enter the decision problem of Equation (\ref{eqn:romlopt}). When $n\rightarrow \infty$, the Rademacher term and the square root terms become zero and the probability of feasibility depends on the asymptotic performance of $\Ialg$ (which converges to $I^*$, the `best-in-class' set function), as desired. The proof is provided in Section \ref{subsec:proof-good-models-mean}.

\section{Robustness guarantee using conditional quantile functions}
\label{sec:quantiles-without-assumption-b}

\begin{theorem} If $\U$ is defined as in Equation (\ref{eqn:empirical-models-quantile-set}), then with probability at least $1 - \delta$ over training sample $S$, we have
\begin{align}
\P_{\{\xUnlabj,\tilde{y}^j\}_{j=1}^{m}}&\Big(F(\pi^*,[\tilde{y}^1...\tilde{y}^m]^T) \in \K\Big) \geq \nonumber\\
 & \left(\left[\frac{1}{n}\sum_{i=1}^{n}\Big(r^{-}_{\epsilon}(y^i - \beta^{\textrm{Alg},\delta_q}(\x^i)) - r^{+}_{\epsilon}(y^i - \beta^{\textrm{Alg},\delta_p}(\x^i))\Big) - \frac{8}{\epsilon}\RA(\B_0) - 2\sqrt{\frac{\log\frac{2}{\delta}}{2n}}\right]_{+}\right)^m,
\end{align}
where $\epsilon > 0$ is a pre-determined constant, $\Big[\cdot \Big]_{+}$ is shorthand for $\max(0,\cdot)$, $r^{-}_{\epsilon}(z) := \min\Big(1,\max\left(0, - \frac{z}{\epsilon}\right)\Big)$ and 
$ r^{+}_{\epsilon}(z) := \min\Big(1,\max\left(0,1 - \frac{z}{\epsilon}\right)\Big)$.
\label{theorem:empirical-models-quantile}
\end{theorem}

The robustness guarantee is established by replacing $\P_{\x,y}(y \leq \beta(\x))$ with the expectation of a related indicator random variable. By majorizing this random variable by random variables defined using functions $r^{-}_{\epsilon}$ and $r^{+}_{\epsilon}$, we were able to use the machinery of Rademacher concentration results. The proof is provided in Section \ref{subsec:proof-empricial-function-quantile}.

\section{Robustness guarantee using a single set of good models}
\label{sec:good-models-mean}

Here we consider the third method prescribed in Section \ref{subsec:good-models-based}. Let $\betaa \in \B_0$ be the model output by the empirical risk minimization procedure. Let empirical risk $l_{S}(\beta) = \frac{1}{n}\sum_{i=1}^n l(\beta(\x^i),y^i)$ be a function of our sample $S$. Let $A$ produce $\betaa$ according to $\betaa \in \arg\min_{\beta \in \B_0}l_{S}(\beta)$. That is, the algorithm $A$ is minimizing the empirical loss. We define $\B$ using the empirical Rademacher average, defined in Section \ref{sec:empirical-functions-general}, as follows: 
\begin{align}
\B := \Bigg\{\beta \in \B_0: l_{S}(\beta)  \leq l_{S}(\betaa) + \left. 2\RA_{S}(l\circ \B_0) + 4M\sqrt{\frac{\log\frac{3}{\delta}}{2n}}\right\},
\label{eqn:learningset}
\end{align}
where $M$ is a bound on the range of the loss function $l$, and $\delta$ is pre-specified and parameterizes the probabilistic guarantee on the robust optimal solution. $\RA_{S}(l\circ \B_0)$ is the empirical Rademacher average of the function class $l\circ\B_0 := \left\{l(\beta(\cdot),\cdot): \beta \in \B_0\right\}$.

We define $\U_{-\B} := E^m$ ($m$ copies of $E$) where $E$ satisfies Equation (\ref{eqn:assumptionA}) for a given $\delta_e$ and $m$ is the number of predictions (equal to the length of the vector $\u_{\beta}$). Intuitively, $\U_{-\B}$ is capturing the support of prediction errors if we knew the `best-in-class' model $\beta^*$. Recall that these definitions of $\U_{\B}$ and $\U_{-\B}$ lead to the set $\U$ in Equation (\ref{eqn:good-models-mean-set}).

\begin{theorem} If $\;\U$ is defined as in Equation (\ref{eqn:good-models-mean-set}), then the following hold:
\begin{enumerate}
\item With probability at least $1 - \delta$, $\beta^* \in \B$.
\item Robust optimal solution $\policy^*$ of Equation (\ref{eqn:romlopt}) is feasible for $\{(\xt^j,\tilde{y}^j)\}_{j=1}^{m}$ with probability at least $(1 - \delta)(1 - \delta_{e})^m$ over $\{(\xt^j,\tilde{y}^j)\}_{j=1}^{m}$ and $S$. That is,
\begin{align*}
\P_{S,\{(\xt^j,\tilde{y}^j)\}_{j=1}^m}\left(F(\policy^*,[\tilde{y}^1 \;\hdots\; \tilde{y}^m]^T) \in \K\right) \geq (1 - \delta)(1 - \delta_{e})^m.
\end{align*}
\end{enumerate}
\label{theorem:good-models-mean}
\end{theorem}
The above theorem holds for any bounded loss function $l$.  It guarantees that $\pi^*$ will be robust to parameter $\u$ with components $ \u_{\beta} = [\beta^*(\xt^1)\hdots \beta^*(\xt^j) \hdots \beta^*(\xt^m)]^T $
and $ \u_{-\beta} = [\tilde{y}^1 \hdots \tilde{y}^j \hdots \tilde{y}^m]^T - [\beta^*(\xt^1)\hdots \beta^*(\xt^j) \hdots \beta^*(\xt^m)]^T$ because the sum of these components is equal to $[\tilde{y}^1 \;\hdots\; \tilde{y}^m]^T$.

We insure against most possible realizations of $\{\tilde{y}^j\}_{j=1}^{m}$  in a particular way: by first ensuring $\beta^*$ belongs to $\B$ with high probability (see Theorem \ref{theorem:good-models-mean} part (1)) and then ensuring that the random errors $\tilde{y}^j- \beta^*(\xt^j)$ are in $\U_{-\B}$ also with high probability. Thus the $\{\tilde{y}^j\}_{j=1}^{m}$ belong to the set $\U_{\B}+ \U_{-\B}$ with high probability. 

This theorem also tells us how the choice of $\B_0$ affects the size of our uncertainty set precursor $\B$. Interestingly enough, if we work with a (possibly infinite) set of predictive models $\B_0$ such that its empirical Rademacher average $\RA_{S}(l\circ \B_0)$ scales as $O(n^{-\frac{1}{2}})$, then we have similar quantitative dependence on $n$ compared to that of confidence-interval based approaches (that make explicit distributional assumptions - see Section \ref{subsec:dependent} - whereas we do not need to make such assumptions). In fact, for many well studied model classes the scaling of the empirical Rademacher average is indeed $O(n^{-\frac{1}{2}})$ which we will review shortly.

One of the advantages of defining uncertainty set precursor $\B$ in the way we proposed is that it directly links the uncertainty in decision making to the loss function $l(\beta(\x),y)$ and sample $S$ of the machine learning step. One advantage of using the empirical Rademacher average in defining $\B$ is that it makes use of the data sample $S$ in its definition, and can reflect the properties of the particular unknown distribution $\P_{\x,y}$ of the data source.


\subsection{Robustness guarantees when the hypothesis set $\B_0$ is finite:}
\label{sec:finite}

When $\B_0$ consists of a finite number of models, we can define $\B$ without using the notion of Rademacher averages. Let $|\B_0|$ represent the size of the set $\B_0$. Then we can define the set of good models as:
\begin{align}
\B := \Bigg\{\beta \in \B_0: l_{S}(\beta)   \leq l_{S}(\betaa) + M\sqrt{\frac{\log|\B_0| + \log\frac{2}{\delta}}{2n}}  + M\sqrt{\frac{\log\frac{2}{\delta}}{2n}}\Bigg\},
\label{eqn:learningsetfinite}
\end{align}
where $n, \delta,M,l_{S}(\cdot)$ and $\betaa$ have the same definitions as before.
\begin{theorem} For finite $\B_0$, the conclusion of Theorem \ref{theorem:good-models-mean} holds if $\;\U$ in Equation (\ref{eqn:good-models-mean-set}) is defined using $\B$ described in Equation (\ref{eqn:learningsetfinite}).
\label{theorem:good-models-mean-finite}
\end{theorem}

\subsection{Constructing $\U$ using PAC-Bayes theory:}
\label{sec:pacbayes}

If the learning step is a classification task, we can also define $\B$ using the PAC-Bayes framework of \citet{mcallester1999pac}, where PAC means ``probably approximately correct''. This framework does not seek a single empirically good classifier $\betaa$ and instead finds a good ``posterior'' distribution $Q$ over the hypothesis set $\B_0$. The corresponding theory provides a probabilistic guarantee on the performance of the classifiers that holds uniformly over all posterior distributions within a class of distributions. The framework then picks a $Q$ using data sample $S$ so that a $Q$-weighted deterministic classifier (or a $Q$-based randomized classifier) has the optimal probabilistic guarantee.

Consider the Q-based (randomized) Gibbs classifier $G_{Q}$, which makes each prediction by choosing a classifier from $\B_0$ according to $Q$. Let the Q-based Gibbs classifier have the following definitions of risks: (a) expected risk $R(G_{Q}) := \E_{\beta \in Q}[l_{\P}(\beta)]$, and (b) empirical risk $R_{S}(G_{Q}) := \E_{\beta \in Q}[l_{S}(\beta)]$ where $l_{\P}(\beta)$ and $l_{S}(\beta)$ are the same as in Section \ref{sec:good-models-mean}. The PAC-Bayes framework guarantees that for all $Q$, $R(G_{Q})$ is bounded by $R_{S}(G_{Q})$ and a term which captures the deviation of $Q$ from a pre-specified `prior' distribution $P$  over $\B_0$ as follows: 

\begin{theorem}\citet[][Theorem 2.1]{germain2009pac}: Let $l(\beta(\x),y):=\mathbf{1}[\beta(\x)\neq y]$. For any $\P_{\x,y}$, any $\B_0$, any prior $P$ on $\B_0$, any $\delta \in (0,1]$ and any convex function $\mathcal{D}:[0,1]^2 \rightarrow \R$, we have 
\begin{align}
\P_{S}\Bigg( \forall Q \textrm{ on } \B_0: \mathcal{D}(R_{S}(G_{Q}),R(G_{Q})) \leq \frac{1}{n}\Bigg[ KL(Q||P) + \log\Big(\frac{1}{\delta}\E_{S}\E_{\beta\sim P}&e^{m\mathcal{D}(l_{S}(\beta),l_{\P}(\beta))} \Big) \Bigg]  \Bigg) \nonumber\\
&\geq  1- \delta,
\label{eqn:pacbound}
\end{align}
where $KL(Q||P) := \E_{\beta \sim Q}[\log\frac{Q(\beta)}{P(\beta)}]$.
\end{theorem}
As shown by \cite{germain2009pac}, for a certain choice of the metric $\mathcal{D}$ the above theorem gives a bound on $R(G_{Q})$ that is proportional to $Cn R_{S}(G_{Q}) + KL(Q||P)$ where $C$ is a pre-specified constant. We can minimize this quantity to get an optimal distribution $Q^\textrm{Alg}$ with a closed form expression: $Q^{\textrm{Alg}}(\beta) = \frac{1}{Z}P(\beta)e^{-Cnl_{S}(\beta)}$ where $Z$ is a normalizing constant.
 
The set of good models $\B$, for the model uncertainty set $\U_{\B}$, can be defined by setting $Q^{\textrm{Alg}}$ to be bigger than a threshold, leading to:
\begin{align*}
\B = \left\{ \beta \in \B_0 :  l_{S}(\beta) \leq \frac{\log P(\beta) - \alpha }{nC} \right\},
\end{align*}
where $\alpha > 0 $ is a fixed constant, $P(\beta)$ is the prior probability density of model $\beta$, and $C$ is a constant that appears in the objective when we solve for $Q^{\textrm{Alg}}$. Intuitively, the set $\B$ includes all models such that their empirical error is bounded in a way that considers their scaled log prior density values. By our construction, if $\beta \in \B$, then $Q^{\textrm{Alg}}(\beta)$ is greater than the threshold $\frac{e^{\alpha}}{Z}$. There is no notion of a `best-in-class' model $\beta^*$ in the PAC-Bayes setting and thus we do not have a guarantee similar to Theorem \ref{theorem:good-models-mean}. Nonetheless, $\B$ is data driven and captures those models that have a high posterior density in $\B_0$. $\U_{B}$ and $\U_{-\B}$ are defined using $\B$ and Equation (\ref{eqn:assumptionA}) in the same way as before and used 
to obtain 
$\pi^*$.


\subsection{Contrasting this method with that of \citet{gi03}:}
\label{subsec:dependent}
\citet{gi03} assume distributional properties on $\hist$ (\textbf{Assumptions GI1}) in addition to assuming a functional form for $\beta(\x)$ (\textbf{Assumption GI2}) while working with robust portfolio selection problems. In particular, let $y = \beta(\x) + \epsilon$, where $\beta(\x) = \betavec^T\x$ is the functional form of the model. Let us assume that $\x^i \in \X \subseteq \R^d$ are chosen by the experimenter and are not random. The only source of randomness is through $\epsilon$ which is independent from example to example and is assumed to be distributed according to $\Nr(0,\sigma^2)$ with variance $\sigma^2$ known.  Then an estimator of $\betavec^*$ (the `best-in-class' model)  is given by:
$\betavec^{\textrm{Alg}} = (X^TX)^{-1}X^TY$,
where $X$ is a matrix with $n$ rows, one for each $\x^i$ and $Y$ is an $n\times1$ vector with the i$^{th}$ element being $y^i$. Here assume that $X^TX$ is invertible. Substituting $Y = X\betavec^* + \boldsymbol{\epsilon}$ in the expression for $\betavec^{\textrm{Alg}}$ gives us: 
$\betavec^{\textrm{Alg}} - \betavec^* = (X^TX)^{-1}X^T\boldsymbol{\epsilon}$,
which is then distributed as $\Nr(0,\sigma^2(X^TX)^{-1})$. Thus, the real-valued function $g(\betavec^*,S) := \frac{1}{\sigma^2}(\betavec^{\textrm{Alg}} - \betavec^*)^{T}(X^TX)(\betavec^{\textrm{Alg}} - \betavec^*)$ is a $\chi^2_d$ distributed random variable. Because of this, we can find a range such that with high probability the $\chi^2_d$ distributed random variable $g(\betavec^*,S)$ belongs to it. We can adapt this approach to our notation by choosing $\B$ based on this interval, giving us an ellipsoid centered at $\betavec^{\textrm{Alg}}$:  $\B = \left\{\betavec: \frac{1}{\sigma^2}(\betavec^{\textrm{Alg}} - \betavec)(X^TX)(\betavec^{\textrm{Alg}} - \betavec) \leq c \right\}$, where $c$ is a constant that determines how much of the probability mass of $\chi^2_d$ is within the set $\B$. 

Set $\U_{-\B}$ can be defined using our assumption about the model residuals: $\epsilon = (y - \betavec^Tx) \sim \Nr(0,\sigma^2)$. In particular, using Equation (\ref{eqn:assumptionA}), we can obtain interval $E = [-e,e]$ for any desired value of $\delta_e$ by solving the equation: $\int_{-e}^{e}\frac{1}{\sqrt{2\pi}\sigma}e^{-\frac{s}{2\sigma^2}}ds = 1 - \delta_e$. 

Using $\B$ (equivalently, $\U_{\B}$) and $\U_{-\B}$ as defined above in the robust problem of Equation (\ref{eqn:romlopt}) gives us a guarantee on the robustness of $\pi^*$ to future realizations of $y$ if \textbf{Assumptions GI1} and \textbf{Assumption GI2} hold. If noise variance $\sigma^2$ is unknown, regression theory provides the following fix: we obtain an unbiased estimator of $\sigma^2$ given by $s^2 = \frac{\|Y-X\betavec^{\textrm{Alg}}\|_2^2}{n-d}$. The resulting scaled random variable $ \frac{1}{ds^2}(\betavec^{\textrm{Alg}} - \betavec)(X^TX)(\betavec^{\textrm{Alg}} - \betavec)$ is $F$-distributed with $d$ degrees of freedom in the numerator and $n-d$ degrees of freedom in the denominator \citep{anderson58}. A set of good models $\B$ can be defined in the same way as before. The constant $c$ now determines how much of the probability mass of an $F_{d,n-d}$-distributed random variable is within $\B$.

Note that both \textbf{Assumptions GI1} and \textbf{Assumption GI2} (or their variations for similar models) are heavily needed to justify these constructions. Contrast this with the setting of Section \ref{sec:good-models-mean} where much weaker assumptions were made and the setting of Section \ref{sec:empirical-functions-general}, where the only assumption made is that the data are drawn i.i.d from some distribution. 
Because our assumptions are much weaker, our result applies to many different loss functions and lends itself naturally to many different machine learning approaches.

\noindent\textbf{Evaluating the empirical Rademacher average}: In the expression for $\B$ in Equation (\ref{eqn:learningset}), it may sometimes be difficult to compute the value of $\RA_{S}(l\circ\B_0)$ efficiently. In these cases, we have two options. The first one involves finding upper bounds on $\RA_{S}(l\circ\B_0)$. This can be tricky as $\RA_{S}$ depends on the data. The second one involves defining $\B$ directly in terms of the \textit{Rademacher average} $\RA(l\circ\B_0)$:
\begin{align}
\B :=& \Bigg\{\beta \in \B_0:
l_{S}(
\beta) \leq l_{S}(\betaa) + 2\RA(l\circ \B_0) + 3M\sqrt{\frac{\log\frac{2}{\delta}}{2n}}\Bigg\}.
\label{eqn:learningset2}
\end{align}
It can be shown that the optimal robust solution obtained using the set in (\ref{eqn:learningset2}) enjoys a guarantee similar to the solution obtained using the set in  (\ref{eqn:learningset}) with different constants. We can make use of the various relationships in Theorem 12 of \cite{bartlett02} to upper bound $\RA(l\circ\B_0)$ or $\RA_{S}(l\circ\B_0)$ analytically. The following are some examples:
\begin{itemize}
\item For linear function classes with squared loss as the loss function, we have:
$
\RA(\B_0) \leq \frac{X_bB_b}{\sqrt{n}} \textrm{, and } 
\RA(l\circ\B_0) \leq 8X_bB_b\frac{X_bB_b}{\sqrt{n}}.
$
where the latter inequality uses Corollary 3.17 in \cite{talagrand91} that relates $\RA(l\circ\B_0)$ and $\RA(\B_0)$. That is, when the loss function $l(\beta(\x),y)$ is $\L$-Lipschitz we have: $\RA(l\circ\B_0) \leq 2\L \cdot \RA(\B_0)$. For the squared loss function, $\L = 4X_bB_b$ if $\forall \x  \in \X, \|\x\|_2 \leq X_b$ and $\forall \beta  \in \B_0, \|\beta\|_2 \leq B_b$. Note that this bound does not depend on data sample $S$.
\item For kernel based function classes with Lipschitz loss functions, $\B_0$ can be written as:
\begin{align*}
\B_0  = \left\{x\mapsto \sum_{i=1}^{n}\alpha^i k(x,x^i):
 n \in \N, x \in \X, \sum_{i,j}\alpha^i\alpha^j K(x^i,x^j) \leq B_b\right\},
\end{align*}
where $k:\X\times\X\rightarrow \R$ is a bounded kernel ($k$ is called a kernel if an $n\times n$ Gram matrix $K$ with entries $(K)_{i,j} = k(x^i,x^j)$ is positive semi-definite). This function class is used in Support Vector Machines (SVMs) \citep[e.g., see][]{shawe-taylor00}
where the loss function is the hinge-loss. The following bound \citep[see][Lemma 22]{bartlett02} applies when the loss function is $\L$-Lipschitz (as is the hinge loss):
\begin{align*}
&\RA_{S}(\B_0) \leq \frac{B_b}{n}\sqrt{\sum_{i=1}^{n}k(x^i,x^j)},\;\;\textrm{ and }\;\; 
\RA_{S}(l\circ\B_0) \leq 2\L\frac{B_b}{n}\sqrt{\sum_{i=1}^{n}k(x^i,x^j)}.
\end{align*}

This upper bound reduces to the previous case (linear function class and squared loss) when we choose the appropriate kernel and loss function. In particular, using the dot product kernel $k(x^i,x^j) = (x^i)^Tx^j$ we get:
\begin{align*}
\RA_{S}(l\circ\B_0) 
\leq 2\L\frac{B_b}{n}\sqrt{\sum_{i=1}^{n}k(x^i,x^j)} = 2\L \frac{B_b}{n}\sqrt{\sum_{i=1}^{n}(x^i)^Tx^j}
\leq 2\L \frac{B_b}{n}\sqrt{nX_b^2}
= 8X_bB_b\frac{X_bB_b}{\sqrt{n}}.
\end{align*}

\end{itemize}

\section{Robustness guarantee using sets of good conditional quantile models}
\label{sec:good-models-quantiles}

Consider the fourth method prescribed in Section \ref{subsec:good-models-based}. Let us define $\B^{\delta_p}$ and $\B^{\delta_q}$ using $\hist$ in a very similar way to defining $\B$ in Equation (\ref{eqn:learningset}). Let the empirical risk minimization procedure using the pinball loss output a conditional quantile model $\betaat$, given sample $S = \hist$ of size $n$ and parameter $\tau$. That is, let $l_{S}^{\tau}(\beta) = \frac{1}{n}\sum_{i=1}^n l^{\tau}(\beta(\x^i),y^i)$ and  $\betaat \in \arg\min_{\beta \in \B_0}l_{S}^{\tau}(\beta)$. The following definition of $\B^{\tau}$ gives us the two sets when $\tau = \delta_p$ and $\tau = \delta_q$:
\begin{align}
\B^{\tau} := 
\left\{\beta 
\in \B_0:  l_{S}^{\tau}(\beta) \leq l_{S}^{\tau}(\betaat) 
  +
  2\RA_{S}(l^{\tau}\circ \B_0) + 4M\sqrt{\frac{\log\frac{3}{\delta}}{2n}}\right\},
\label{eqn:quantileset}
\end{align}
where $M$ is a bound on the range of the loss function $l^{\tau}$, $\delta$ is a pre-specified constant and $\RA_{S}(l^{\tau}\circ \B_0)$ is the empirical Rademacher average of the function class $l^{\tau}\circ\B_0 := \{\beta\mapsto l^{\tau}(\beta(\cdot),\cdot): \beta \in \B_0\}$.
The guarantee on the robust optimal solution of Equation (\ref{eqn:romlopt}) is given by the following theorem.

\begin{theorem}
If $\U$ is defined as described in Equation (\ref{eqn:good-models-quantile-set}), using sets $\B^{\delta_p}, \B^{\delta_q}$ defined in Equation (\ref{eqn:quantileset}) and set $E$ in Equation (\ref{eqn:assumptionB}) along with \textbf{Assumption B}, then the following hold:
\begin{enumerate}
\item With probability at least $1 - \delta$, $\beta^{\tau,*} \in \B^{\tau}$ for $\tau = \delta_p$ and $\tau = \delta_q$ individually.
\item Robust optimal solution $\policy^*$ of Equation (\ref{eqn:romlopt}) is feasible for $\{(\xt^j,\tilde{y}^j)\}_{j=1}^{m}$ with probability at least $(1 - \delta)\left[(1 - \delta_e^{\delta_p})^m + (1 - \delta_e^{\delta_q})^m\right] + (\delta_q - \delta_p)^m - 2$ over $\{(\xt^j,\tilde{y}^j)\}_{j=1}^{m}$ and $S$. That is,
\begin{align*}
\P_{S,\{(\xt^j,\tilde{y}^j)\}_{j=1}^m}
\left(F(\policy^*,[\tilde{y}^1 ... \tilde{y}^m]^T) \in \K\right) \geq
(1 - \delta)\left[(1 - \delta_e^{\delta_p})^m + (1 - \delta_e^{\delta_q})^m\right] + (\delta_q - \delta_p)^m - 2.
\end{align*}
\end{enumerate}
\label{theorem:good-models-quantile}
\end{theorem}

In the theorem, the guarantee follows from designing $\U$ such  that the predictions made by the `best-in-class' conditional quantile functions $\beta^{\delta_p,*},\beta^{\delta_q,*}$ and their residuals are captured in each interval defining $\U$. This ensures that the realization  $[\tilde{y}^1 ... \tilde{y}^m]^T \in \U$ with high probability.

The guarantees in Sections \ref{sec:good-models-mean} and this section do not assume anything about the form of the source distribution. These bounds do what learning theory is designed to do \citep{bousquet03}, which is provide insight into the important quantities for learning and how they scale. More importantly, here they provide insight beyond prediction, specifically into robustness for decision making. We generally do not use learning theoretic bounds directly in practice (e.g., SVMs do not minimize the generalization bounds that motivated their derivation). To translate our results in practice, we suggest using our workflow to construct the uncertainty sets as in Equations (\ref{eqn:learningset}), (\ref{eqn:learningsetfinite}) or (\ref{eqn:learningset2}), replacing the Rademacher average term with an appropriate choice of parameters. A practitioner can also perform a type of sensitivity analysis for our approach by varying the size of the uncertainty sets and assessing the corresponding results.

\subsection{Insights and Comparison of Main Results}

Before we move onto the proofs, we recap the main results. Theorem \ref{theorem:empirical-models-generic} provides a very general results that pertains to any algorithm. Intuitively, it states that as long as the algorithm's result is robust to most of the training examples, and as long as the algorithm can only produce simple functions, it will likely be robust to all points in the test set. This is true for any unknown distribution of data, with no assumptions on the distribution.

Theorem \ref{theorem:empirical-models-generic} does not provide any insights on how to construct an algorithm for data-driven robust optimization, since it holds for any algorithm. Theorem \ref{theorem:empirical-models-quantile}, on the other hand, provides a result that holds for quantile regression methods. Here we use an algorithm that produces an estimate for a lower quantile and an estimate for a higher quantile, and chooses the policy to be robust to all points between these quantile estimates. The result applies to any method for producing such estimates. It states that, for this choice of policy, the solution will be robust to all points on the test set with high probability. The bound will be tighter when the class of quantile estimation functions produced by the algorithm is simpler. Theorem \ref{theorem:empirical-models-quantile} is close to being a special case of Theorem \ref{theorem:empirical-models-generic}. Theorem \ref{theorem:empirical-models-quantile}'s loss function is similar to a special case of Theorem \ref{theorem:empirical-models-generic}'s, and the complexity term differs only through a Lipschitz constant of the loss function which is explicitly taken into account in Theorem \ref{theorem:empirical-models-quantile} but not in Theorem \ref{theorem:empirical-models-generic}.

Theorems \ref{theorem:good-models-mean} and \ref{theorem:good-models-quantile} rely on mild probabilistic assumptions that can make the bounds tighter. These theorems consider the full set of ``good" models, that is, models with small loss on the training set, and expand outwards to include more points into the uncertainty set; thus these theorems take into account both the behavior on the training set and the assumed behavior on the full distribution of data.

For the assumption underlying Theorems \ref{theorem:good-models-mean} and \ref{theorem:good-models-quantile}, recall that $\delta_e$ is the probability that the tails of the distribution are within $E$ of the true mean or quantile estimates. There is a tradeoff in assumptions between $E$ and $\delta_e$, in the sense that the policy needs to be robust to a larger uncertainty set if $E$ is large; larger $E$ leads to conservative policy choices. At the same time, if $E$ is larger, our assumption that $E$ includes the tails of the distribution should be stronger, leading to smaller $\delta_e$. When $\delta_e$ is smaller, the probabilistic guarantee on robustness is also stronger.

Theorem \ref{theorem:good-models-mean}'s result holds for any algorithm that produces estimates of centrality for $y$ given $x$ (e.g., mean or median). Theorem \ref{theorem:good-models-quantile}'s result holds for any algorithm that produces quantile estimates. We believe that the uncertainty set construction used for Theorem \ref{theorem:good-models-quantile} is the most natural ones to use, regardless of whether the assumptions relating $\delta_e$ and $E$ hold precisely. To recap, this is where we compute the highest estimate of the upper quantile from all good models, compute the lowest estimate of the lower quantile from all good models, and expand outwards, to produce the uncertainty set.

\section{Proofs}
\label{sec:proofs}

Before we proceed with the proofs of guarantees for the four methods in Sections \ref{sec:empirical-functions-general}-\ref{sec:good-models-quantiles}, we state an intermediate result we will make use of in all four proofs. This result gives a uniform  probabilistic guarantee on the deviation between empirical loss and expected loss of prediction models in terms of the Rademacher average. It holds for any set of models $\mathcal{F}$ and a bounded
loss function $l$.

\begin{lemma} With probability at least $1-\delta$ over sample $S$,  
\begin{align*}
\max_{f \in \mathcal{F}} \left|l_{\P}(f) - l_{S}(f)\right| \leq 2\RA(l\circ\mathcal{F}) + M\sqrt{\frac{\log\frac{1}{\delta}}{2n}}.
\end{align*}
\label{lemma:unifdeviation}
\end{lemma}
\textit{Proof.} Here, $\max_{f \in \mathcal{F}} |l_{\P}(f) - l_{S}(f)|$ is a random variable that depends on the sample $S$ through $l_{S}()$. We can use the (one-sided) McDiarmid's inequality to claim that this random variable is close to its mean as $n$ increases. 
\begin{lemma} \emph{McDiarmid's inequality} \citep{mcdiarmid1989method}:
Let $z^1,...,z^n$ be $n$ i.i.d. random variables in a set A and $h(z^1,...,z^n)$ be a function such that for all $i=1,...,n$
\begin{align*}
\sup_{(z^1,...,z^n,\tilde{z}) \in A^{n+1}}|h(z^1,...,
z^i,...,z^{n}) 
- h(z^1,...,\tilde{z},...,z^{n})| \leq c.
\end{align*}
\begin{align*}
\textrm{Then for all }\epsilon >0, \;\; \P_{z^1,...,z^n}\Bigg(h(z^1,...,z^n) -& \E[h(z^1,...,z^n)] > \epsilon \Bigg)
\leq 
\exp\left(-\frac{2\epsilon^2}{nc^2}\right).
\end{align*}
\label{lemma:mcdiarmid}
\end{lemma}
In our case, the function $h$ is $\max_{f \in \mathcal{F}} \left|l_{\P}(f) - l_{S}(f)\right|$. We can show 
that if the $i^{th}$ instance in the sample $S$ is perturbed, the maximum change in the function value is $\frac{M}{n}$: We first consider the case when  $\max_{f \in \mathcal{F}} \left|l_{\P}(f) - l_{S}(f)\right| \geq \max_{f \in \mathcal{F}} \left|l_{\P}(f) - l_{S^i}(f)\right|$. Here $l_{S^i}(f)$ is the same as $l_{S}(f)$ except for the $i^{\textrm{th}}$ example, which is changed from $(\x^i,y^i)$ to a new example $\x_{\circ}^{i},y_{\circ}^i$. Also let $f^{\circ} \in \arg\max_{f \in \mathcal{F}} \left|l_{\P}(f) - l_{S}(f)\right|$. Then,
\begin{align*}
\max_{f \in \mathcal{F}} &\left|l_{\P}(f) - l_{S}(f)\right| - \max_{f \in \mathcal{F}} \left|l_{\P}(f) - l_{S^i}(f)\right| \\
\leq&  \left|l_{\P}(f^{\circ}) - l_{S}(f^{\circ})\right| - \left|l_{\P}(f^{\circ}) - l_{S^i}(f^{\circ})\right| \;\;\;\textrm{(because $f^\circ$ may not maximize the second term)}\\
\leq&   \left| - l_{S}(f^{\circ})  + l_{S^i}(f^{\circ}) \right| \;\;\;\textrm{(by triangle inequality)}\\
=& \frac{1}{n} \left| l(f^{\circ}(\x_{\circ}^i),y^i)- l(f^{\circ}(\x^i),y^i) \right| 
\leq \frac{M}{n} \;\;\; \textrm{(canceling all except the $i^{\textrm{th}}$ term)}.
\end{align*}
We can do an identical calculation to get the same upper bound $\frac{M}{n}$ if $\max_{f \in \mathcal{F}} \left|l_{\P}(f) - l_{S}(f)\right| \leq \max_{f \in \mathcal{F}} \left|l_{\P}(f) - l_{S^i}(f)\right|$. Thus, with probability at least $1 - \delta$,
\begin{align}
\max_{f \in \mathcal{F}} &\left|l_{\P}(f) - l_{S}(f)\right| \leq 
\E[\max_{f \in \mathcal{F}} \left|l_{\P}(f) - l_{S}(f)\right|] + M\sqrt{\frac{\log\frac{1}{\delta}}{2n}}. \label{eqn:closeness_sub2}
\end{align}

The quantity $\E[\max_{f \in \mathcal{F}} \left|l_{\P}(f) - l_{S}(f)\right|]$ captures the complexity or size of $\mathcal{F}$ (actually, its composition with the loss function $l$, the set $l\circ \mathcal{F}$). We can upper bound this quantity in terms of a Rademacher average 
 using a symmetrization trick.
\begin{lemma} (Upper bound)
\begin{align}
\E[\max_{f \in \mathcal{F}} \left| l_{\P}(f) - l_{S}(f)\right|] \leq 2\RA(l\circ\mathcal{F}). \label{eqn:closeness_sub3}
\end{align}
\end{lemma}
\textit{Proof.} See Theorem 8 in \cite{bartlett02} for essentially a similar claim.

Substituting for $\E[\max_{\beta \in \B_0} \left|l_{\P}(\beta) - l_{S}(\beta)\right|]$ from (\ref{eqn:closeness_sub3}) into (\ref{eqn:closeness_sub2}) gives us the desired result.\qed

\subsection{Proof of Theorem \ref{theorem:empirical-models-generic}}
\label{subsec:proof-empricial-models-generic}

According to Lemma \ref{lemma:unifdeviation}, the following holds with probability at least $1-\delta$ over sample $S$,  
\begin{align*}
\max_{f \in \mathcal{F}} \left|l_{\P}(f) - l_{S}(f)\right| \leq 2\RA(l\circ\mathcal{F}) + M\sqrt{\frac{\log\frac{1}{\delta}}{2n}}.
\end{align*}
where $\mathcal{F}$ is a set of models, $l$ is a bounded loss function (bounded by $M$), $l_\P(f) = \E_{\x,y}[l(f(x),y)]$, $l_{S}(f) = \frac{1}{n}\sum_{i=1}^{n}l(f(\x^i),y^i)$ and $\RA(l\circ\mathcal{F}) = \E_{S,\sigma}[\sup_{f \in \mathcal{F}}\frac{1}{n}|\sum_{i=1}^{n}\sigma^il(f(\x^i),y^i) |]$.

We can apply this lemma to the case when $\mathcal{F} = \I$ (that is, $f(x) = I(x)$) and $l(I(x),y) = \ind[y \notin I(x)]$. The range of the loss function is $[0,1]$, which is a bounded set. Thus, with probability at least $1-\delta$ over sample $S$,
\begin{align*}
&\max_{I \in \I} \left|l_{\P}(I) - l_{S}(I)\right| \leq 2\RA(l\circ\I) + \sqrt{\frac{\log\frac{1}{\delta}}{2n}},\\
\textrm{or equivalently,}&\\
&\forall I \in \I : l_{S}(I) - 2\RA(l\circ\I) - \sqrt{\frac{\log\frac{1}{\delta}}{2n}} \leq l_{\P}(I) \leq l_{S}(I) + 2\RA(l\circ\I) + \sqrt{\frac{\log\frac{1}{\delta}}{2n}}.
\end{align*}
The above is a uniform convergence statement. Since it holds for $\Ialg$ as well, we can state the following: with probability at least $1-\delta$ over $S$,
\begin{align}
l_{S}(\Ialg) - c(\delta) \leq  \P_{\x,y}(y \notin \Ialg(\x)) \leq l_{S}(\Ialg) + c(\delta), \textrm{ or equivalently,} \nonumber\\
1 - \left(l_{S}(\Ialg) - c(\delta)\right) \geq  \P_{\x,y}(y \in \Ialg(\x)) \geq 1- \left(l_{S}(\Ialg) + c(\delta)\right),
\label{eqn:general-set-prob-ineq}
\end{align} 
where $c(\delta) =  2\RA(l\circ\I) + \sqrt{\frac{\log\frac{1}{\delta}}{2n}}$ and we use the relation $l_{\P}(I) = \E_{\x,y}[\ind[y \notin I(\x)]] = \P_{\x,y}(y \notin I(\x))$. 

The second inequality in Equation (\ref{eqn:general-set-prob-ineq}) gives a lower bound on the probability that an unseen label belongs to the interval specified by the function $\Ialg(\x)$. We can extend this lower bound to $m$  unseen new realizations $\{\tilde{y}^j\}_{j=1}^{m}$ as follows. With probability at least $1-\delta$ over $S$,
\begin{align*}
\P_{\xUnlabj,\tilde{y}^j}(\tilde{y}^j \in \Ialg(\xUnlabj)) \geq 1- \left(l_{S}(\Ialg) + c(\delta)\right); \; j=1,...,m.
\end{align*}
Then, with probability $\geq 1-\delta$ over $S$,
\begin{align*}
\P_{\{\xUnlabj,\tilde{y}^j\}_{j=1}^{m}}([\tilde{y}^1,...,\tilde{y}^m]^T \in \Pi_{j=1}^{m}\Ialg(\xUnlabj)) \geq (1- \left(l_{S}(\Ialg) + c(\delta)\right))^{m},
\end{align*}
where we used the fact that these $m$ events $\{\tilde{y}^j \in \Ialg(\xUnlabj)\}, j=1,...,m$ are mutually independent given sample $S$. 

Note that if $[\tilde{y}^1,...,\tilde{y}^m]^T \in \Pi_{j=1}^{m}\Ialg(\xUnlabj)$, then the robust optimal solution $\policy^*$ is feasible for the future label realizations $\{\tilde{y}^j\}_{j=1}^{m}$ because it is feasible for each of the $m$ elements in $\U = \Pi_{j=1}^{m}\Ialg(\xUnlabj)$ by definition. This gives us the desired feasibility result on $\policy^*$. \qed

\subsection{Proof of Theorem \ref{theorem:empirical-models-quantile}}
\label{subsec:proof-empricial-function-quantile}

As in the previous proof, according to Lemma \ref{lemma:unifdeviation} the following holds with probability at least $1-\delta$ over sample $S$:  
\begin{align*}
\max_{f \in \mathcal{F}} \left|l_{\P}(f) - l_{S}(f)\right| \leq 2\RA(l\circ\mathcal{F}) + M\sqrt{\frac{\log\frac{1}{\delta}}{2n}}.
\end{align*}
where $\mathcal{F}$ is a set of models, $l$ is a bounded loss function (bounded by $M$), $l_\P(f) = \E_{\x,y}[l(f(x),y)]$, $l_{S}(f) = \frac{1}{n}\sum_{i=1}^{n}l(f(\x^i),y^i)$ and $$\RA(l\circ\mathcal{F}) = \E_{S,\sigma}\left[\sup_{f \in \mathcal{F}}\frac{1}{n}\left|\sum_{i=1}^{n}\sigma^il(f(\x^i),y^i) \right|\right].$$

We will apply the lemma in two cases. For both cases, let the model set $\B_0$ be the set of conditional quantile models. For the first case, let the loss function be $r^{-}_{\epsilon}(y - \beta(\x))$ and for the second case, let the loss function be $r^{+}_{\epsilon}(y - \beta(\x))$ (both functions are defined in the statement of Theorem \ref{theorem:empirical-models-quantile}). The range of both functions is $[0,1]$, and thus bounded.
Further, since Lemma \ref{lemma:unifdeviation} is a uniform deviation statement, the inequality also holds for model $\beta^{\textrm{Alg},\tau}$, derived from sample $S$ (say, by minimizing the pinball loss), with probability at least $1- \delta $. Thus, we have the following two probabilistic statements:
\begin{itemize}
\item With prob. $\geq 1-\delta$ over $S$,
\begin{align}
\left|\E_{\x,y}[r^{-}_{\epsilon}(y -  \beta^{\textrm{Alg},\tau}(\x))] -\frac{1}{n}\sum_{i=1}^{n}r^{-}_{\epsilon}(y^i - \beta^{\textrm{Alg},\tau}(\x^i))\right| \leq 2\RA(r^{-}_{\epsilon}\circ\B_0) + \sqrt{\frac{\log\frac{1}{\delta_2}}{2n}}. \label{eqn:concrminus}
\end{align}
\item With prob. $\geq 1-\delta$ over $S$,
\begin{align}
\left|\E_{\x,y}[r^{+}_{\epsilon}(y -  \beta^{\textrm{Alg},\tau}(\x))] -\frac{1}{n}\sum_{i=1}^{n}r^{+}_{\epsilon}(y^i - \beta^{\textrm{Alg},\tau}(\x^i))\right| \leq 2\RA(r^{+}_{\epsilon}\circ\B_0) + \sqrt{\frac{\log\frac{1}{\delta_2}}{2n}}. \label{eqn:concrplus}
\end{align}
\end{itemize}
From these inequalities, we get the following lemma \citep[similar to][Theorem 7]{takeuchi2006nonparametric}:

\begin{lemma}
With probability at least $1-\delta$ over sample $S$, the following inequalities hold separately:
\begin{align}
\frac{1}{n}\sum_{i=1}^{n}r^{-}_{\epsilon}(y^i - \beta^{\textrm{Alg},\tau}(\x^i)) - c \leq &\P_{\x,y}(y \leq \beta^{\textrm{Alg},\tau}(\x)),\textrm{ and} \label{eqn:ineqlower}\\
&\P_{\x,y}(y \leq \beta^{\textrm{Alg},\tau}(\x))  \leq \frac{1}{n}\sum_{i=1}^{n}r^{+}_{\epsilon}(y^i - \beta^{\textrm{Alg},\tau}(\x^i)) + c, \label{eqn:inequpper}
\end{align}
where $c := \frac{4}{\epsilon}\RA(\B_0) + \sqrt{\frac{\log\frac{1}{\delta}}{2n}}$.
\label{lemma:quantileineq}

\end{lemma}
\textit{Proof} (of Lemma \ref{lemma:quantileineq})
From the Ledoux-Talagrand contraction inequality, we know $\RA(l\circ \B_0) \leq 2\L\RA(\B_0)$. In our case, both $r^{-}_{\epsilon}$ and $r^{+}_{\epsilon}$ have Lipschitz constant equal to $1/\epsilon$. Let $c$ be defined as in the statement of the lemma. From the inequalities (\ref{eqn:concrminus}) and (\ref{eqn:concrplus}) we get the one sided inequalities:
\begin{align*}
\E_{\x,y}[r^{-}_{\epsilon}(y -  \beta^{\textrm{Alg},\tau}(\x))] \geq \frac{1}{n}\sum_{i=1}^{n}r^{-}_{\epsilon}(y^i - \beta^{\textrm{Alg},\tau}(\x^i)) - c, \textrm{ and}\\
\E_{\x,y}[r^{+}_{\epsilon}(y -  \beta^{\textrm{Alg},\tau}(\x))] \leq \frac{1}{n}\sum_{i=1}^{n}r^{+}_{\epsilon}(y^i - \beta^{\textrm{Alg},\tau}(\x^i)) + c.
\end{align*}

Further, for any $\beta$, we can bound $\P_{\x,y}(y \leq \beta(\x)) = \E_{\x,y}[\ind[y \leq \beta(\x)]] $ from both sides because of the following inequalities:
\begin{align}
\E_{\x,y}[\ind[y \leq \beta(\x)]] &\leq \E_{\x,y}[r^{+}_{\epsilon}(y - \beta(\x))], \textrm{ and} \label{eqn:majorize}\\
\E_{\x,y}[\ind[y \leq \beta(\x)]] &\geq \E_{\x,y}[r^{-}_{\epsilon}(y - \beta(\x))]. \label{eqn:minorize}
\end{align}

Thus we get:
\begin{itemize}
\item with prob. $\geq 1-\delta$, $\P_{\x,y}(y \leq \beta^{\textrm{Alg},\tau}(\x)) \geq \frac{1}{n}\sum_{i=1}^{n}r^{-}_{\epsilon}(y^i - \beta^{\textrm{Alg},\tau}(\x^i)) - c$, and
\item with prob. $\geq 1-\delta$, $\P_{\x,y}(y \leq \beta^{\textrm{Alg},\tau}(\x)) \leq \frac{1}{n}\sum_{i=1}^{n}r^{+}_{\epsilon}(y^i - \beta^{\textrm{Alg},\tau}(\x^i)) + c$.
\end{itemize}
\qed

Continuing with the proof of Theorem \ref{theorem:empirical-models-quantile}, we
 apply Lemma \ref{lemma:quantileineq} with $\tau = \delta_p$ within inequality (\ref{eqn:inequpper}) and with $\tau = \delta_q$ within inequality (\ref{eqn:ineqlower}), where $1 \leq \delta_p < \delta_q \leq 1$ to obtain:
\begin{itemize}
\item with prob. $\geq 1-\delta$, $\P_{\x,y}(y \leq \beta^{\textrm{Alg},\delta_q}(\x)) \geq \frac{1}{n}\sum_{i=1}^{n}r^{-}_{\epsilon}(y^i - \beta^{\textrm{Alg},\delta_q}(\x^i)) - c$, and
\item with prob. $\geq 1-\delta$, $\P_{\x,y}(y \leq \beta^{\textrm{Alg},\delta_p}(\x)) \leq \frac{1}{n}\sum_{i=1}^{n}r^{+}_{\epsilon}(y^i - \beta^{\textrm{Alg},\delta_p}(\x^i)) + c$.
\end{itemize}
The bounds hold with probabilities $1-\delta$ each implying that they together hold with probability $1-2\delta$. Now, 
\begin{align*}
\P_{\x,y}(&\beta^{\textrm{Alg},\delta_p}(\x) < y \leq \beta^{\textrm{Alg},\delta_q}(\x))\\
&= \P_{\x,y}(\{\beta^{\textrm{Alg},\delta_p}(\x) < y \} \cap \{ y \leq \beta^{\textrm{Alg},\delta_q}(\x)\})\\
&= 1 - \P_{\x,y}(\{y \leq \beta^{\textrm{Alg},\delta_p}(\x) \} \cup \{\beta^{\textrm{Alg},\delta_q}(\x) < y \})\\
&\geq 1 - \left(\P_{\x,y}(y \leq \beta^{\textrm{Alg},\delta_p}(\x)) + \P_{\x,y}(\beta^{\textrm{Alg},\delta_q}(\x) < y )\right)\\
&= 1 - \left(\P_{\x,y}(y \leq \beta^{\textrm{Alg},\delta_p}(\x)) + 1 - \P_{\x,y}(y \leq \beta^{\textrm{Alg},\delta_q}(\x))\right)\\
&=  \P_{\x,y}(y \leq \beta^{\textrm{Alg},\delta_q}(\x)) - \P_{\x,y}(y \leq \beta^{\textrm{Alg},\delta_p}(\x))\\
&\overset{(*)}{\geq} \frac{1}{n}\sum_{i=1}^{n}r^{-}_{\epsilon}(y^i - \beta^{\textrm{Alg},\delta_q}(\x^i)) - \frac{1}{n}\sum_{i=1}^{n}r^{+}_{\epsilon}(y^i - \beta^{\textrm{Alg},\delta_p}(\x^i)) - 2c,
\end{align*}
where in step $(*)$, we substituted upper and lower bounds of the two random variables of $S$, $\P_{\x,y}(y \leq \beta^{\textrm{Alg},\delta_q}(\x))$ and $\P_{\x,y}(y \leq \beta^{\textrm{Alg},\delta_p}(\x))$. Thus, with probability $\geq 1-2\delta$ over $S$,
\begin{align*}
\P_{\x,y}(y \in [\beta^{\textrm{Alg},\delta_p}(\x) ,\beta^{\textrm{Alg},\delta_q}(\x)]) \geq \frac{1}{n}\sum_{i=1}^{n}\left(r^{-}_{\epsilon}(y^i - \beta^{\textrm{Alg},\delta_q}(\x^i)) - r^{+}_{\epsilon}(y^i - \beta^{\textrm{Alg},\delta_p}(\x^i))\right) - 2c.
\end{align*}
In the above statement, we have a lower bound on the probability that a new unseen realization $y$ belongs to the random interval $[\beta^{\textrm{Alg},\delta_p}(\x) ,\beta^{\textrm{Alg},\delta_q}(\x)]$.

We can extend this lower bound to the setting of $m$ simultaneous lower bounds corresponding to $m$ unseen new realizations $\{\tilde{y}^j\}_{j=1}^{m}$ in our decision problem as follows. We know that with probability $\geq 1-2\delta$ over $S$,
\begin{align*}
\P_{\xUnlabj,\tilde{y}^j}(\tilde{y}^j \in [\beta^{\textrm{Alg},\delta_p}(\xUnlabj) ,\beta^{\textrm{Alg},\delta_q}(\xUnlabj)]) \geq \Delta(S); \; j=1,...,m,
\end{align*}
where $\Delta(S) := \frac{1}{n}\sum_{i=1}^{n}\left(r^{-}_{\epsilon}(y^i - \beta^{\textrm{Alg},\delta_q}(\x^i)) - r^{+}_{\epsilon}(y^i - \beta^{\textrm{Alg},\delta_p}(\x^i))\right)  - 2c$.
Then, with probability $\geq 1-2\delta$ with respect to sample $S$,
\begin{align*}
\P_{\{\xUnlabj,\tilde{y}^j\}_{j=1}^{m}}([\tilde{y}^1,...,\tilde{y}^m]^T \in \Pi_{j=1}^{m}[\beta^{\textrm{Alg},\delta_p}(\xUnlabj) ,\beta^{\textrm{Alg},\delta_q}(\xUnlabj)]) \geq \Delta(S)^m,
\end{align*}
where we used the fact that these $m$ events $\{\tilde{y}^j \in [\beta^{\textrm{Alg},\delta_p}(\xUnlabj) ,\beta^{\textrm{Alg},\delta_q}(\xUnlabj)]\}, j=1,...,m$ are mutually independent given sample $S$.

Note that if $[\tilde{y}^1,...,\tilde{y}^m]^T \in \Pi_{j=1}^{m}[\beta^{\textrm{Alg},\delta_p}(\xUnlabj) ,\beta^{\textrm{Alg},\delta_q}(\xUnlabj)]$, then it also belongs to $\U$ defined by Equation (\ref{eqn:empirical-models-quantile-set}). Further, the robust optimal solution $\policy^*$ will be feasible for $\{\tilde{y}^j\}_{j=1}^{m}$ because it is feasible for every element in $\U$ by definition. Thus, changing $\delta$ to $\delta/2$ (with an appropriate change in the constant $c$ in Equations (\ref{eqn:ineqlower}) and (\ref{eqn:inequpper})) gives us the desired feasibility result on $\policy^*$. \qed

\subsection{Proof of Theorem \ref{theorem:good-models-mean}}
\label{subsec:proof-good-models-mean}

Consider the term $l_{S}(\beta^*) - l_{S}(\betaa)$, which depends on the random sample $S$. We can upper bound it by:
\begin{align}
l_{S}(\beta^*) - &l_{S}(\betaa) \nonumber\\
			& = l_{S}(\beta^*) - l_{\P}(\beta^*) + l_{\P}(\beta^*) - l_{S}(\betaa)\nonumber\\
			& \leq l_{S}(\beta^*) - l_{\P}(\beta^*) + l_{\P}(\betaa) - l_{S}(\betaa)\nonumber\\
                         & \leq l_{S}(\beta^*) - l_{\P}(\beta^*) + \max_{\beta \in \B_0} \left|l_{\P}(\beta) - l_{S}(\beta)\right| \label{eqn:closeness}
\end{align}
where we added and subtracted $ l_{\P}(\beta^*)$ in the first step, then in the second step substituted $\betaa$ for $\beta^*$ in the third term to increase the value of the right hand side, and finally in the last step, replaced the last two terms with an absolute max operation over $\B_0$.

The first term in the expression on the right hand side of (\ref{eqn:closeness}) will go to zero in probability as $n \rightarrow \infty$ due to concentration, and this can be quantified for finite $n$ via Hoeffding's inequality.
\begin{lemma} \emph{(One-sided Hoeffding's inequality.)}
Let $z^1,...,z^n$ and $z$ be i.i.d. random variables and let $h$ be a bounded function, $a\leq h(z) \leq b$. Then for all $\epsilon > 0$ we have
\begin{align*}
\P_{z^1,...,z^n}\left( \frac{1}{n}\sum_{i=1}^{n} h(z^i) - \right.&\left. \E_{z}[h(z)] > \epsilon \right) 
\leq \exp\left(-\frac{2n\epsilon^2}{(b-a)^2}\right).
\end{align*}
\label{lemma:hoeffding}
\end{lemma}
In our case, the sample $S$ is represented by $\hist$. The function $l(\beta^*(x),y)$ is bounded in the interval $[0,M]$. Thus the empirical mean $\frac{1}{n}\sum_{i=1}^{n} l(\beta^*(x^i),y^i)$ ($= l_S(\beta^*)$) gets close to its mean $\E[l(\beta^*(x),y)] $ ($= l_{\P}(\beta^*)$) as $n$ increases. In particular, we see that with probability at least $1-\delta_1$,
\begin{align}
l_S(\beta^*) - l_{\P}(\beta^*) \leq M\sqrt{\frac{\log\frac{1}{\delta_1}}{2n}}. \label{eqn:closeness_sub1}
\end{align}

The second term (\ref{eqn:closeness}) can be bounded using Lemma \ref{lemma:unifdeviation} which states that with probability at least $1-\delta$ over sample $S$,  
\begin{align*}
\max_{f \in \mathcal{F}} \left|l_{\P}(f) - l_{S}(f)\right| \leq 2\RA(l\circ\mathcal{F}) + M\sqrt{\frac{\log\frac{1}{\delta}}{2n}}.
\end{align*}
In our case, we set $\mathcal{F} = \B_0$ and $f(x) = \beta(x)$ and $\delta = \delta_2$.

The empirical Rademacher average $\RA_{S}(l\circ\B_0)$ also concentrates around its mean $\RA(l\circ\B_0)$ and this can be proved again by McDiarmid's inequality. In this case, from Lemma \ref{lemma:mcdiarmid}, the function $h$ is represented by $\RA_{S}(l\circ \B_0)$. We can again show  \citep[][Theorem 11]{bartlett02} that if the $i^{th}$ instance in the sample $S$ is perturbed, the maximum change in the function value is $\frac{M}{n}$. 
Thus, with probability at least $1 - \delta_3$,
\begin{align}
\RA(l\circ\B_0) \leq \RA_{S}(l\circ\B_0) + M\sqrt{\frac{\log\frac{1}{\delta_3}}{2n}}. \label{eqn:closeness_sub4}
\end{align}

In summary we have the following statements for the terms on the right hand side of (\ref{eqn:closeness}):
\begin{enumerate}
\item With probability  at least $1-\delta_1$ over $S$, $l_S(\beta^*) - l_{\P}(\beta^*) \leq M\sqrt{\frac{\log\frac{1}{\delta_1}}{2n}}$ from (\ref{eqn:closeness_sub1}).
\item With probability  at least $1-\delta_2$ over $S$, 
\begin{align*}
\max_{\beta \in \B_0} \left|l_{\P}(\beta) - l_{S}(\beta)\right| \leq 2\RA(l\circ\B_0) + M\sqrt{\frac{\log\frac{1}{\delta_2}}{2n}}.
\end{align*} 
\item With probability  at least $1-\delta_3$ over $S$, $\RA(l\circ\B_0) \leq \RA_{S}(l\circ\B_0) + M\sqrt{\frac{\log\frac{1}{\delta_3}}{2n}} $ from (\ref{eqn:closeness_sub4}).
\end{enumerate}

Consider the three corresponding events:
$E_1 = \Bigg\{S \;:\; l_S(\beta^*) - l_{\P}(\beta^*) \leq M\sqrt{\frac{\log\frac{1}{\delta_1}}{2n}}\Bigg\}$, $E_2 = \Bigg\{S \;:\; \max_{\beta \in \B_0} \left(l_{\P}(\beta) - l_{S}(\beta)\right) \leq 2\RA(l\circ\B_0) + M\sqrt{\frac{\log\frac{1}{\delta_2}}{2n}}\Bigg\}$, and
$E_3 = \Bigg\{S \;:\; \RA(l\circ\B_0) \leq \RA_{S}(l\circ\B_0) + M\sqrt{\frac{\log\frac{1}{\delta_3}}{2n}}\Bigg\}$.
We know that with probabilities $\delta_1,\delta_2,\delta_3$ over the random sample $S$, these events do not happen. Thus using the union bound, 
$ \P_{S}(E_1 \cap E_2 \cap E_3) \geq 1 - \delta_1 + \delta_2 + \delta_3
$.
Substituting $\frac{\delta}{3}$ for $\delta_1,\delta_2$ and $\delta_3$ and using (\ref{eqn:closeness}), we that with probability at least $1 - \delta$,
\begin{align*}
l_{S}(\beta^*) - l_{S}(\betaa) \leq  2\RA_{S}(l\circ\B_0) + 4M\sqrt{\frac{\log\frac{3}{\delta}}{2n}}. 
\end{align*}


The implication of this is that the empirical risk for the `best-in-class' function $\beta^*$ is less than the right hand side quantities, all of which are computable. This implies that even though we do not know $\beta^*$, we know it belongs to our uncertainty set precursor $\B$ defined in Equation (\ref{eqn:learningset}) with high probability. In particular, we see that $\beta^* \in \B$ with probability  at least $1 - \delta$ over sample $S$. This is part (1) in the statement of the Theorem.

Part (1) further implies that with probability at least $1-\delta$, $u_{\beta^*} \in \U_{\B}$, and this is true for any $\{\xt^j\}_{j=1}^{m}$. Next we turn our focus toward model residuals. We can extend the probabilistic statement in Equation (\ref{eqn:assumptionA}) to the setting where we have $m$ simultaneous errors using the mutual independence assumption. Thus we have, with probability at least $ (1- \delta_e)^m$ over $\{(\xt^j,\tilde{y}^j)\}_{j=1}^{m}$,
$\max_{j=1,...,m}|\tilde{y}^j-\beta^*(\xt^j)| \in E$.  Using the definition of set $\U_{-\B}$, which is equal to $E^m$, we see that $u_{-\beta*} \in \U_{-\B}$ with probability at least $(1-\delta_e)^m$ over $\{(\xt^j,\tilde{y}^j)\}_{j=1}^{m}$.

We know that the robust optimal solution $\policy^*$ of Equation  (\ref{eqn:romlopt}) is robust to any element of $\U = \U_{\B} + \U_{-\B}$ by definition. In particular, if $\beta^* \in \B$ and $u_{-\beta^*} \in \U_{-\B}$, then $\policy^*$ will be robust to the random vector $\u_{\beta^*} + \u_{-\beta^*}$ (which equals $[\tilde{y}^1\hdots\tilde{y}^m]^T$). 

To get a guarantee of robustness of $\pi^*$ to $\{\tilde{y}^j\}_{j=1}^{m}$, we can combine the two probabilistic statements above (one with respect to $S$ and the other with respect to $\{(\xt^j,\tilde{y}^j)\}_{j=1}^{m}$) using the mutual independence assumption ($S$ and $\{(\xt^j,\tilde{y}^j)\}_{j=1}^{m}$ are mutually independent) as follows:
\begin{align*}
\P_{S,\{(\xt^j,\tilde{y}^j)\}_{j=1}^m}\left(F(\policy^*,[\tilde{y}^1\hdots \tilde{y}^m]^T) \in \K\right) \geq (1-\delta)(1 - \delta_{e})^m.\qed
\end{align*} 

\subsection{Proof of Theorem \ref{theorem:good-models-mean-finite}}
\label{subsec:proof-good-models-mean-finite}

It is sufficient to show that with probability at least $1-\delta$, $\beta^* \in \B$ where $\B$ is defined in Equation (\ref{eqn:learningsetfinite}). To see this, consider the deviation $l_{S}(\beta^*) - l_{S}(\betaa)$. This can be upper bounded in a similar way as in the beginning of the proof of Theorem \ref{theorem:good-models-mean}:
\begin{align*}
l_{S}(\beta^*) - l_{S}(\betaa) \leq l_{S}(\beta^*) - l_{\P}(\beta^*) + \max_{\beta \in \B_0} (l_{\P}(\beta) - l_{S}(\beta)).
\end{align*}
We will upper bound the two deviation terms appearing on the right hand side of the above inequality. Both terms are functions of the random sample $S$.

Lets begin with the term $\max_{\beta \in \B_0} (l_{\P}(\beta) - l_{S}(\beta))$.  We can bound the probability of the event $\{\max_{\beta \in \B_0} (l_{\P}(\beta) - l_{S}(\beta)) > \epsilon\}$ as follows:
\begin{align*}
\P_{S}\Big(&\max_{\beta \in \B_0} (l_{\P}(\beta) - l_{S}(\beta)) > \epsilon\Big)
= \P_{S}\Big(\cup_{i=1}^{|\B_0|}\{l_{\P}(\beta^i) - l_{S}(\beta^i) > \epsilon\}\Big)\\
&\stackrel{(a)}{\leq} \sum_{i=1}^{|\B_0|}\P_{S}\Big(l_{\P}(\beta^i) - l_{S}(\beta^i) > \epsilon\Big) \stackrel{(b)}{=} \sum_{i=1}^{|\B_0|}e^{-\frac{2n\epsilon^2}{M^2}} = e^{\log|\B_0|-\frac{2n\epsilon^2}{M^2}}.
\end{align*}
Here, (a) follows from taking a union bound, and (b) follows from applying Hoeffding's inequality to each fixed model $\beta^i,i=1,...,|\B_0|$.  Setting $\delta_2 = e^{\log|\B_0|-\frac{2n\epsilon^2}{M^2}}$ and replacing $\epsilon$ gives us the following equivalent way to state the same result: with probability at least $1 - \delta_2$ over $S$,
\begin{align*}
\max_{\beta \in \B_0} (l_{\P}(\beta) - l_{S}(\beta)) \leq M\sqrt{\frac{\log|\B_0| + \log(\frac{1}{\delta_2})}{2n}} .
\end{align*}

From Equation (\ref{eqn:closeness_sub1}), we have the following upper bound for the term $l_{S}(\beta^*) - l_{\P}(\beta^*)$ : with probability  at least $1-\delta_1$ over $S$, $l_S(\beta^*) - l_{\P}(\beta^*) \leq M\sqrt{\frac{\log\frac{1}{\delta_1}}{2n}}$. 

Using a union bound with these two observations gives us the following statement when we set $\delta_1 = \delta_2 = \delta/2$: with probability at least $1 - \delta$ over $S$,
$l_{S}(\beta^*) - l_{S}(\betaa) \leq  M\sqrt{\frac{\log|\B_0| + \log(\frac{2}{\delta})}{2n}} + M\sqrt{\frac{\log\frac{2}{\delta}}{2n}}$.
Thus $\beta^* \in \B$ with probability  at least $1-\delta$ as desired.

\qed

\subsection{Proof of Theorem \ref{theorem:good-models-quantile}}
\label{subsec:proof-good-models-quantile}

Proof of part (1) is the same as that of part (1) in Theorem \ref{theorem:good-models-mean}. That is, using the definition of $\B^{\tau}$ in Equation (\ref{eqn:quantileset}) and Lemma \ref{lemma:unifdeviation} with the pinball loss function $l^{\tau}$ we see that $\beta^{\tau,*} \in \B^{\tau}$ with probability  at least  $1-\delta$ over $S$. Thus, part (1) holds when $\tau$ is set to $\delta_p$ and $\delta_q$ individually.

For part (2), we use mutual independence and union bound arguments, similar to part (2) in Theorem \ref{theorem:good-models-mean}. In particular,
\begin{itemize}

\item With prob. $\geq$   $1 - \delta$ over $S$, simultaneously for all $j=1,..,m$, $\beta^{\delta_p,*}(\xUnlabj) \in [\inf\{\beta(\xUnlabj): \beta \in \B^{\delta_p}\},\sup\{\beta(\xUnlabj): \beta \in \B^{\delta_p}\}]$  for any $\{\xt^j\}_{j=1}^{m}$ (from part (1)).

\item With prob. $\geq$   $1 - \delta$ over $S$, simultaneously for all $j=1,..,m$, $\beta^{\delta_q,*}(\xUnlabj) \in [\inf\{\beta(\xUnlabj): \beta \in \B^{\delta_q}\},\sup\{\beta(\xUnlabj): \beta \in \B^{\delta_q}\}]$ for any $\{\xt^j\}_{j=1}^{m}$ (from part (1)).

\item With prob. $\geq$   $(1 - \delta_e^{\delta_p})^m$ over $\{(\xt^j,\tilde{y}^j)\}_{j=1}^{m}$, simultaneously for all $j=1,..,m$, $\mu^{\delta_p}(\xUnlabj) - \beta^{\delta_p,*}(\xUnlabj) \in [-\sup E^{\delta_p}, \sup E^{\delta_p}]$  (using mutual independence assumption and Equation (\ref{eqn:assumptionB})).

\item  With prob. $\geq$   $(1 - \delta_e^{\delta_q})^m$ over $\{(\xt^j,\tilde{y}^j)\}_{j=1}^{m}$, simultaneously for all $j=1,..,m$, $\mu^{\delta_q}(\xUnlabj) - \beta^{\delta_q,*}(\xUnlabj) \in [-\sup E^{\delta_q}, \sup E^{\delta_q}]$ (using mutual independence assumption and Equation (\ref{eqn:assumptionB})).

\end{itemize}

We can again use the mutual independence between $S$ and $\{(\xt^j,\tilde{y}^j)\}_{j=1}^{m}$ to claim the following:

\begin{itemize}

\item With prob. $\geq$ $(1 - \delta)(1 - \delta_e^{\delta_p})^m$ over $S$ and $\{(\xt^j,\tilde{y}^j)\}_{j=1}^{m}$, simultaneously for all $j=1,...,m$, $\mu^{\delta_p}(\xUnlabj) \in [\inf\{\beta(\xUnlabj): \beta \in \B^{\delta_p}\}-\sup E^{\delta_p}, \sup\{\beta(\xUnlabj): \beta \in \B^{\delta_p}\} + \sup E^{\delta_p}]$.
\item With prob. $\geq$ $(1 - \delta)(1 - \delta_e^{\delta_q})^m$ over $S$ and $\{(\xt^j,\tilde{y}^j)\}_{j=1}^{m}$, simultaneously for all $j=1,...,m$, $\mu^{\delta_q}(\xUnlabj) \in [\inf\{\beta(\xUnlabj): \beta \in \B^{\delta_q}\}-\sup E^{\delta_q}, \sup\{\beta(\xUnlabj): \beta \in \B^{\delta_q}\} + \sup E^{\delta_q}]$.
\end{itemize}

We use the general identity from De Morgan's laws and the union bound that if $\P(A_1) \geq c_1$ and $\P(A_2)\geq c_2$, then $\P(A_1\cap A_2) \geq c_1 + c_2 -1$. Applying this to the two events above, we see that with probability at least $
(1 - \delta)\left[(1 - \delta_e^{\delta_p})^m + (1 - \delta_e^{\delta_q})^m\right] -1$ over $S$ and $\{(\xt^j,\tilde{y}^j)\}_{j=1}^{m}$,
\begin{align*}
[\mu^{\delta_p}(\xUnlabj),&\mu^{\delta_q}(\xUnlabj)] \subseteq\\
 [&\inf\{\beta(\xUnlabj): \beta \in \B^{\delta_p}\cup\B^{\delta_q}\}-\sup E^{\delta_p}\cup E^{\delta_q}, \\
 &\sup\{\beta(\xUnlabj): \beta \in \B^{\delta_p}\cup\B^{\delta_q}\}+\sup E^{\delta_p}\cup E^{\delta_q}].
\end{align*}
We also know that simultaneously for all $j$, $\tilde{y}^j$ belongs to $ [\mu^{\delta_p}(\xUnlabj),\mu^{\delta_q}(\xUnlabj)]$ with probability at least $(\delta_q - \delta_p)^m$ over $\{(\xt^j,\tilde{y}^j)\}_{j=1}^{m}$ (mutual independence and definition of conditional quantile function). Thus, again using the identity based on De Morgan's laws and the union bound, we get that with probability at least $(1 - \delta)\left[(1 - \delta_e^{\delta_p})^m + (1 - \delta_e^{\delta_q})^m\right] + (\delta_q - \delta_p)^m - 2$ over $S$ and $\{(\xt^j,\tilde{y}^j)\}_{j=1}^{m}$, $[\tilde{y}^1 ...\tilde{y}^m]^T $ belongs to the set
\begin{align*}
\Pi_{j=1}^{m}\Big[&\inf\{\beta(\xUnlabj): \beta \in \B^{\delta_p}\cup\B^{\delta_q}\}-\sup E^{\delta_p}\cup E^{\delta_q}, \\
 &\sup\{\beta(\xUnlabj): \beta \in \B^{\delta_p}\cup\B^{\delta_q}\}+\sup E^{\delta_p}\cup E^{\delta_q}\Big].
\end{align*}
Since $\U$ is defined precisely using the above product set, 
we conclude that the robust optimal solution $\policy^*$ is feasible for $\{\tilde{y}^j\}_{j=1}^{m}$ with the desired guarantee. \qed

\section{Conclusion}
\label{sec:conclusion}
In this work, we presented two principled approaches (four methods) of constructing uncertainty sets for robust optimization based on statistical learning theory. These methods can be used broadly for data-driven robust optimization, and apply to any problem where the data are drawn from an unknown distribution. The first two methods can be applied without any distributional assumptions, and the other two methods require very mild distributional assumptions, which is that the user knows one statistic about the tail of the distribution. The results in this paper show that statistical learning theory, derived for guarantees on prediction quality of statistical models, can be used for guarantees on the robustness of an optimization problem.\\

\ACKNOWLEDGMENT{Funding for this project comes in part from Ford-MIT Alliance and NSF grant IIS-1053407.}

\bibliographystyle{plainnat}
\bibliography{roml_paper}
\end{document}